\input amssym.def 
\input epsf.tex 
\input eplain
\font\piccolo=cmr10 at 8truept
\font \piccobo=cmbx10 at 8truept
 at 12truept

\def \hb{\hfill \break}  
\def \pf{\noindent{$\underline{\hbox{Proof}}$.\ }}

\def \muu{\overline {\mu}}

\def \wdd{\widetilde {D}}
\def \wss{\widetilde {S}}

\def \wii{\widetilde {I}}
 
\def \wxx{\widetilde {X}}
\def \wff{\widetilde {f}}
\def \fh#1 {\Bbb F_{#1} }  
\def \cvd{\hfill {$\diamond$} \smallskip}
\def \grc {\hbox{\bf \piccobo GRC}}

\centerline{\bf On the singularities of surfaces ruled by conics}
\medskip
\centerline{ Michela Brundu\footnote{$^*$}{\piccolo Corresponding author. E-mail: brundu@units.it}, Gianni Sacchiero}

\centerline{\piccolo Dipartimento di Matematica e Geoscienze, Universit\`a di Trieste, Via Valerio 12/1 - 34127 Trieste, Italy}

\bigskip
\centerline{\piccobo Abstract}

{\piccolo We classify the singularities of a surface ruled by conics:  they are rational double points of type $A_n$ or $D_n$. This is proved by showing that they arise from a precise series of blow--ups of a suitable surface geometrically ruled by conics. We determine also the family of such surfaces which are birational models of a given surface ruled by conics and obtained in a ``minimal way" from it.}

\medskip
{\piccobo Keywords: } {\piccolo  projective surface, conic bundle, ruled surface, rational singularities.}

\smallskip
{\piccolo {\piccobo MSC:}  14 J 26, 14 E 05, 14 D 06}

\medskip
\noindent
{\bf Introduction}

Projective surfaces ruled by conics arise naturally in the study of the moduli
space of four--gonal curves
${\cal M}_{g,4} \subset {\cal M}_g$ inside of the moduli space of curves of genus $g$. Indeed, it is known that
the canonical model of such curves lies on a relative
hyperquadric $S$ in a three-dimensional rational normal scroll $V=
\Bbb P ({ \cal E})$,
that is a divisor of the type $2H - \beta F$, with $H$ tautological divisor and  $F$ a fiber.

\smallskip

In analogy with the Maroni invariant in the trigonal case, 
 the splitting type  of  the vector bundle ${\cal E}= {\cal O}(a) \oplus {\cal O}(b)  \oplus {\cal O}(g-3-a-b)$ gives  a first  description of the curve. But in [3], using the invariant $t$ of a precise birational  smooth model ${\Bbb F}_t$ of $S$ and a further invariant $\lambda$ (which is substantially the minimum degree of a linear series of the curve, out of the four--gonal one), a more precise description is given, essentially by  linking  the properties of the surface $S$ with the geometry of the curve (for instance the authors prove that $\deg(S) =g+\lambda -t -5$). More precisely, with a suitable use of these four invariant, they define and describe a stratification in irreducible locally closed subsets ${\cal M}_g^{{\lambda},t}(a,b)$
of ${\cal M}_{g,4}$, which has been extensively studied in the cited paper.

\smallskip

In this frame, becomes essential a precise description of the surface $S$, its singular locus and the type of its singularities.
This motivates the local and global study of projective surfaces with isolated singularities,
fibered over ${\Bbb P}^1$, with general fiber a smooth conic.

\smallskip

We dealt with a very similar subject in [2], but the investigation carried out in the present paper  differs from and is more general than the previous article in two important respects. First of all, the
previous study considered only rational surfaces, whereas here we
treat the general case. Second, it is known that every surface $S$
ruled by conics is birational to a surface which is geometrically
ruled by conics and everywhere smooth. This means
that $S$ can be obtained from such a surface by means of a finite
number of blow-ups and blow-downs (see [6], Ch.~4, Sect.~3). In the
previous study, we restricted to the case of surfaces which arise in
this way using blow-ups at {\it distinct} points and the main theorem (1.9) of [2]  classifies the
singularities of a surface ruled by conics in that situation. 

\smallskip
In the present paper, the main result (Theorem 2.4) describes {\it all} the possible singularities of a surface ruled by conics.

\smallskip
Let us point out that the nature of the  classification given in the present paper is
of an algorithmic type, and in principle could be used for computational purposes.

\smallskip 

Surfaces ruled by conics are roughly regarded as a special case of conic bundles. This subject has been widely developed in the literature, but mainly in the case of  the dimension of the base variety is at least two (only recall the well--known papers [10],  [7], [8], among the several works spread in decades). In Section 0 we first remind the basic notions and the general results concerning conic bundles contained in the paper of Sarkisov [9] which apply also when the base of the conic bundle has dimension one. Then we
prove that a surface ruled by conics can be
regarded as an embedded conic bundle, at least if it has isolated singularities.

\smallskip

In Sect.1 we describe in detail how the singularities of a surface ruled by conics can be
computed {\it via} a series of blow-ups and contractions from a
surface {\it geometrically} ruled by conics. 
In Sect.2  we give a complete classification
of {\it all $ $}   surfaces ruled by conics, by showing that the singularities computed in the previous section are the all possible ones.
\smallskip
We also solve the inverse problem: given a surface ruled
by conics, how can one recover, in a minimal way, a
birationally--equivalent surface which is geometrically ruled by
conics? We determine the (finite) family of such birational models
in Sect.~3. This last result will play an important role in the study of the moduli space of four--gonal curves cited before.

\bigskip
\noindent
{\bf 0. Preliminary notions and known results} 

\medskip
All varieties in this paper are assumed to be algebraic over a fixed algebraically closed field of
characteristic $0$. By {\it projective surface} we mean an irreducible and reduced surface in $\Bbb P^N$.

\medskip
\noindent
{\bf Notation.} 
\item{$(a)$} If $X$ is an algebraic
surface and $P \in X$ is any point, then $\sigma_P$ will denote the {\it blow--up
morphism} at the point $P$ and the obtained surface, the {\it blow--up} of $X$ at $P$,
will be denoted by $Bl_P(X)$ (or $\wxx$ if the centre of the monoidal transformation is
clear). Briefly:
$$
\sigma_P : \; Bl_P(X) = \wxx \longrightarrow X.
$$ 
\item{$(b)$} If $\wxx$ is is an algebraic
surface and $E \subset \wxx$ is a contractible curve, we will
set:
$$
con(E) : \; \wxx \longrightarrow X
$$
the {\it blow--down morphism} giving the contraction of $E$ to a point of $X$.
\item{$(c)$} If $X$ is as before and $D$ is a divisor on $X$, then
$$
\Phi_D: \quad X \longrightarrow X' \subset \Bbb P(H^0({\cal O}_X(D)))
$$
denotes the morphism associated to $D$.

\bigskip
\noindent
{\bf Definition 0.1.}
Let $C$ be a smooth irreducible curve of genus $g$ and let $\cal F$ be a rank 2 vector bundle on $C$. If
$D$ is a very ample $k$--secant divisor on $\Bbb P ({\cal F})$, then the surface 
$S_0:= \Phi_D(\Bbb P ({\cal F}))$ is said to be {\it geometrically  $k$--ruled} over $C$.
Equivalently,  a projective surface
$S_0 \subset \Bbb P^N$ is  geometrically $k$--ruled  over  $C$ if there exists
a surjective morphism $\pi: \; S_0 \longrightarrow C$
such that the fibre $\pi^{-1}(y)$ is a smooth rational curve of degree $k$ for every point $y \in C$.

\medskip
In particular, a geometrically  $k$--ruled  surface $S_0$ is smooth hence $\pi: \; S_0 \longrightarrow C$ is flat, being $C$ a smooth curve (see [1],  pg. 91). 

\bigskip
\noindent
{\bf Definition 0.2.} Let $C$ be as before. We say that  a projective surface
$S \subset \Bbb P^N$ is   {\it $k$--ruled} over  $C$ if there exists a surjective morphism 
$\pi: \; S \longrightarrow C$ and  a non--empty open subset
$U\subseteq C$ such that: 
\item \item {-} $\pi$ is flat with  fibres  $\pi^{-1}(y)$ of degree $k$ and arithmetic genus
$0$ for every point $y \in C$; 
\item \item {-} the fibre $\pi^{-1}(y)$ is  smooth  for every point $y \in U$.
\par

\bigskip
In this paper we will treat $2$--ruled surfaces, hence the   very
ample divisor $D$ on $\Bbb P({\cal F})$ is bisecant. 
In this case the surface
 $S_0=\Phi_D(\Bbb P({\cal F}))$ will be called {\it geometrically ruled by conics} while the surface $S$  will be said {\it ruled by conics}. \hb
In Sections 1 and 2 we show that a surface $S$ ruled by conics is birational to a surface $S_0$ geometrically ruled by conics and in Section 3 we analyze these birational models of $S$.

\bigskip
It is clear that a surface $S$ ruled by conics corresponds to a curve in $\Bbb P^5$, the space parametrizing the conics of the plane. Since this curve can be locally approximated to a line, then the surface $S$ can be locally expressed as a pencil of conics. In this way one can see that if $P$ is a singular point  of $S$ then necessarily three facts occur: $P$ is a base--point of the above pencil of conics, the fibre $F_P$ of $S$ containing $P$ is singular at $P$ and, finally, $F_P$ contains at most one more singular point of $S$. \hb
This leads to the following basic fact:

\medskip
\proclaim
Proposition 0.3.  Let $\pi: S \rightarrow C$ be a surface ruled by conics on a smooth curve $C$. Then $S$ has only isolated singularities.
\par
\pf
Let $\Sigma$ be a one--dimensional component of the singular locus $S_{sing}$. If $\Sigma$ intersect all fibres of $S$, this implies  (from the above observation) that each fibre is a singular conic. But this is impossible. \hb
Hence $\Sigma$ contains a fibre; so such fibre is contained in the singular locus of $S$. But, again from the previous remark, each fibre contains at most two singular points of $S$. So also this case cannot occur. \hb
Therefore $S_{sing}$ is a zero--dimensional subset of $S$.
\cvd

\bigskip
There is an extensive literature on conic bundles; here we mention some notions and results collecting them mainly from the work of Sarkisov (see [9]).

\bigskip
\noindent
{\bf Definition 0.4.} Let $S$ and $C$ be irreducible algebraic varieties and $C$ be non singular. A triple $(S,C,\pi)$, where $\pi: \; S \longrightarrow C$  is a rational map whose generic fibre is an irreducible rational curve, is called a {\it conic bundle over the base $C$}.

\bigskip
\noindent
{\bf Definition 0.5.} A conic bundle $(S,C,\pi)$ is called {\it regular} if  $\pi$ is a flat morphism of nonsingular varieties. 

\bigskip
Let $(S,C,\pi)$ be a conic bundle. If $y \in C$, denote by $S_y$ the fibre $\pi^{-1}(y)$ and, if  $\cal F$ is a sheaf on $S$, denote by ${\cal F}_y$ the restriction of $\cal F$  to the fibre $S_y$. Finally, set $\omega_S$ and $K_S$, respectively, the canonical sheaf and the canonical divisor of $S$.

\bigskip
\noindent
{\bf Remark 0.6.} Let $(S,C,\pi)$ be a regular conic bundle. Then ${\cal O}_S(-K_S)$ is flat over $C$, since for all $x \in S$, $\left( {\cal O}_S(-K_S)  \right)_x \cong {\cal O}_{x,S}$, which is flat over ${\cal O}_{\pi(x),C}$ by assumption.

\bigskip
\noindent
{\bf Remark 0.7.}  Let $(S,C,\pi)$ be a regular conic bundle. Then $H^i(S_y, {\cal O}_S(-K_S)_y) = 0$ for $i \ge 1$, for each  $y \in C$. 
To show this, note that the cohomology groups vanish for all $i \ge 2$ since $\dim(S_y)=1$. \hb
Let ${\cal N}_{y/C}$ denote the normal bundle at the point $y \in C$; then, since $C$ is smooth, ${\cal N}_{y/C}$ is free 
free of rank $r:=\dim(C)$. On the other hand,  $\pi^{*}({\cal N}_{y/C}) = {\cal N}_{S_y/S}$, so also the latter is (locally) free of  rank $r$. Therefore
$$
\bigwedge ^{r} \,{\cal N}_{S_y/S} = {\cal O}_{S_y}.
$$
Since $S_y \subset S$ is a nonsingular subvariety of codimension $r $, it holds (by [6], Ch. II, 8.20 and the above equality) that
$$
\omega_{S_y} \cong \omega_S \otimes \bigwedge ^{r} \,{\cal N}_{S_y/S} = \omega_S \otimes {\cal O}_{S_y}.
$$
Finally, dualizing the above equality, we obtain that the restriction to the fibre of the anticanonical sheaf is exactly the anticanonical sheaf of the fibre:
$$
{\cal O}_S(-K_S) \otimes  {\cal O}_{S_y} \cong {\cal O}_{S_y}(-K_{S_y}) . 
$$
But  each fibre of $\pi$ is isomorphic to $\Bbb P^1$, hence 
$$
H^1(S_y, {\cal O}_S(-K_S)_y) = H^1(S_y,{\cal O}_{S_y}(-K_{S_y})) = H^1(\Bbb P^1,{\cal O}_{\Bbb P^1}(2)) = 0.
$$

\medskip
\proclaim
Proposition 0.8.   If  $(S,C,\pi)$ is a regular conic bundle, then  we have:
\item{- } $R^i \pi_*({\cal O}_S(-K_S)) = 0$ for all $i \ge 1$;
\item{- } $R^0 \pi_*({\cal O}_S(-K_S))$ is a locally free sheaf of rank 3.
\par
\pf
Using  0.6 and 0.7, we can apply [6], Ex 11.8, Ch. III  , obtaining that $R^i \pi_*({\cal O}_S(-K_S)) = 0$ for all $i \ge 1$ (in a neighborhood of each $y \in C$). \hb
With the same argument, it is easy to see that $\dim H^0 (S_y, {\cal O}_S(-K_S)_y)$ is constant  on $C$. So, using Grauert Theorem ([6], 12.9, Ch. III), we obtain that 
$R^0 \pi_*({\cal O}_S(-K_S))$ is locally free on $S$ and 
$$
(R^0 \pi_*({\cal O}_S(-K_S)))_y \cong H^0 (S_y, {\cal O}_S(-K_S)_y) \cong H^0(\Bbb P^1,{\cal O}_{\Bbb P^1}(2))
$$
which is a vector space of dimension $3$.
\cvd

\goodbreak

\medskip
Setting ${\cal E} := R^0 \pi_*({\cal O}_S(-K_S))$ the above locally free sheaf of rank $3$ on $C$ and taking into account the above isomorphism ${\cal E}_y \cong H^0 (S_y, {\cal O}_S(-K_S)_y)$, it is clear that the morphism associated to the divisor $-K_S$ of $S$ is an embedding $S \hookrightarrow \Bbb P({\cal E})$ under which the image of each fibre $S_y$ of $\pi$ is a conic in the corresponding fibre $\Bbb P({\cal E})_y \cong \Bbb P^2$. For this reason, $S$ is said {\it conic} bundle.\hb
The above construction, performed in the regular conic bundle case, leads to the described embedding $S \hookrightarrow \Bbb P({\cal E})$. But one can have an analogous embedding also in a more general case:

\medskip
\noindent
{\bf Definition 0.9.}  Let  ${\cal E}$ be a locally free sheaf of rank $3$ on $C$ and $\tau: \; \Bbb P({\cal E}) \longrightarrow C$ be the natural projection. An irreducible reduced divisor $S \subset  \Bbb P({\cal E}) $ such that the triple $(S, C, \tau_{|S})$ is a conic bundle is called an {\it embedded conic bundle}.

\medskip
The case of $1$--dimensional base is not particularly focused in [9], but it turns out that in this frame can be settled also the surfaces ruled by conics (in the sense of 0.1 and 0.2), as the following result shows:

\medskip
\proclaim
Proposition 0.10.  Let $S \subset \Bbb P^N$ be a surface ruled by conics over the curve $C$ via $\pi: S \rightarrow C$. Then:
\item{i)} each singular point of $S$ is a Gorenstein singularity;
\item{ii)}  $S$ is normal;
\item{iii)} $(S, C, \pi)$ is an embedded conic bundle.
\par
\pf $i)$ For each $P \in S$, let $T_P(S)$ be the tangent space to $S$ at $P$, $F_P$ the fibre of $S$ containing 
$P$ and $U_P$ a unisecant on $S$ passing through $P$. Then
$$
\dim(T_P(S)) \le \dim (T_P(F_P))  +  \dim (T_P(U_P)).
$$
Since $C$ is a smooth curve, then $\dim (T_P(U_P)) =1$, while $\dim (T_P(F_P))$ is $1$ or $2$ accordingly if $P$ is a smooth or singular point of the fibre $F_P$. This is due to the fact that $F_P$ has degree $2$ and $F_P$ has no embedded points since $\pi$ is flat.
Therefore $\dim(T_P(S)) \le 3$; in particular, if  $S$ is not smooth at $P$ then $P$ is a hypersurface singularity of $S$, hence a Gorenstein singularity. \hb
$ii)$ It is immediate from a well--known theorem (see for instance [11]): if $P$ is a two--dimensional singularities, then $P$ is normal if and only if it is Cohen--Macaulay and an isolated singularity. \hb
$iii)$ It is a consequence of [4], Proposition 2.1:  if $\pi: \; X \longrightarrow Z$ is a morphism, $X$ is a Gorenstein scheme, $Z$ is smooth and the fibres of $\pi$ are all (possibly degenerate) conics in $\Bbb P^2$, then $X$  is a flat conic bundle. 
\cvd

\medskip
The last notion we want to recall about conic bundles is that of {\it degeneration divisor}. In [9] it is defined as the vanishing locus of a homomorphism of certain cohomology groups (see [9], 1.6). Even if the author gives a complete characterization of such divisor in the case of regular conic bundle (in particular when $S$ is smooth), let us recall the properties which hold also in the general case: the degeneration divisor $\Delta$ of a conic bundle $(S,C,\pi)$ is such that, for all points $y \in C \setminus \Delta$, the conic $\pi^{-1}(y)$ has rank 3, i.e. it is smooth. If $y \in \Delta$ then the corresponding conic has rank $\le 2$. Finally, if $\pi$ is flat,  the rank of each conic is non--zero.

\medskip
In a forthcoming paper of the moduli space of $4$--gonal curves, we need a very detailed description of all the possible singularities of $S$ and the relations between the invariant of its non--singular models and the invariants of $\Bbb P({\cal E})$. As far as we know, there are not such results in the literature. Hence, in what follows, we determine
a procedure (a sort of  constructive algorithm) in order to detect and classify the singularities of a projective surface ruled by conics.

\bigskip
\noindent
{\bf 1. Singularities arising from elementary transformations of ``main type" } 

\medskip
Since, as observed in the previous Section, any $2$--ruled surface $S$ has a finite number of singular points, then it can be obtained by a suitable geometrically
$2$--ruled surface $S_0= \Phi_D(\Bbb P ({\cal F}))$ by a finite number of monoidal
transformations. In other words, denoting by
$\wss_0$ the surface obtained by a sequence $\sigma$ of blow--ups of $S_0$, 
setting  $\wdd$ to be the strict transform of $D$ via $\sigma$ and 
$B$  the base locus of $\wdd$, then $S$ can be obtained in the following way:
$$
\commdiag{ 
\wss_0 &  \mapright^\sigma & S_0 \cr
\mapdown \lft{\Phi_{\wdd -B}} \cr
S
}
\eqno {(1)}
$$
where $\Phi_{\wdd-B}$ is a birational morphism.

\goodbreak

Note that, for each rank $2$ vector
bundle $\cal F$, the scroll $\Bbb P({\cal F})$ can be locally expressed as $U \times \Bbb P^1$, where $U$ is
an open affine subset of the base curve
$C$. Moreover, also the choice of the very ample $2$--secant divisor $D$ on $\Bbb P({\cal F})$ does not
affect $S_0$ since $\Phi_D$ is an isomorphism.

\medskip
Let us explicitly describe $\sigma$ as composition of a chain of blow--ups centered in suitable points:  we assume that the centre consists of points either belonging to the same fibre or infinitely near to it. \hb 
 Consider a point  
$P_1 \in S_0$ and let $f_0:= \pi^{-1}(y)$ be the fibre of $S_0$ containing $P_1$. Let
us consider the blow--up of $S_0$ at $P_1$ and the corresponding projection on $C$, say
$\pi_1$: 
$$
     \harrowlength=50pt
     \sarrowlength=20pt
     \commdiag{
     Bl_{P_1} (S_0) :=   S_1  &\mapright^{\sigma_{P_1}}& S_0 \supset \pi^{-1}(y) = f_0 \ni P_1\cr
    &\arrow(1,-1)\lft{\pi_1} \quad  \arrow(-1,-1)\rt{\pi}\cr
     &\qquad C\ni y \cr 
     }
 $$
Denote also by $f_1 := \pi_1^{-1}(y)$ the total transform of $f_0$ via $\sigma_{P_1}$.
\hb 
Take now $P_2 \in f_1$ and consider the corresponding blow--up $\sigma_{P_2}: \; S_2
\longrightarrow S_1$. With obvious notations, we can repeat this construction and 
obtain a sequence of blow--ups:
$$
 \harrowlength=40pt
\commdiag{
\wss_0:= S_n & \mapright^{\sigma_{P_n}} & S_{n-1} & \harrowlength=20pt \mapright& \cdots & \harrowlength=20pt \mapright& S_2 & \mapright^{\sigma_{P_2}} &  S_1  &  \mapright^{\sigma_{P_1}} & S_0 \cr
\hfill \cup \; & & & & & &\cup && \cup \; && \cup \;\cr
\wff_0:=f_n & & & & & & f_2 & \hfill P_2 \in & f_1 &  \hfill  P_1 \in & f_0 \cr}
\eqno(2)
$$
where, for all $i$, we set $P_i \in f_{i-1}$, $f_i:= \pi_i^{-1}(y)$ and 
$\pi_i: S_i:= Bl_{P_i}(S_{i-1}) \longrightarrow C$  the natural projection.

\medskip
\noindent
{\bf Remark 1.1.} The surface $\wss_0$ is   smooth at
each point of
$f_n$. Moreover, since $\sigma:= \sigma_{P_n} \circ \cdots \circ \sigma_{P_1}$ is an
isomorphism from
$\wss_0
\setminus f_n$ to
$S_0
\setminus f_0$, then $\wss_0$ is smooth everywhere.

\medskip
\noindent
{\bf Definition 1.2.} Let $\wss_0 \supset f_n$ be as in $(2)$.
Then we say that $f_n$ is a fibre of {\it
level $n$ over $f_0$}.

\medskip
Let us consider the single fibres of the surfaces involved in diagram $(1)$:
$$
\commdiag{ 
\wss_0 &  \mapright^\sigma & S_0 \cr
\mapdown \lft{\Phi_{\wdd -B}} \cr
S
}
\qquad
\qquad
\qquad
\commdiag{ 
f_n &  \mapright^\sigma & f_0 \cr
\mapdown \lft{\Phi_{\wdd -B}} \cr
F
}
\eqno {(3)}
$$

Clearly,  the fibre $F$ of $S$ is uniquely determined by $f_n$.

\medskip
\noindent
{\bf Definition 1.3.} We say that the fibre $F \subset S$
is an {\it embedded fibre of level
$n$} if
$$
n = \min_i \; \{ \hbox{there exists a blow--up $\sigma: \wss_0 \rightarrow S_0$ and a fibre 
$f_i \subset \wss_0$  of level
$i$ such that} \; F =
\Phi_{\wdd -B}(f_i)
\}.
$$

The purpose of this section is to describe some of the possible fibres of $\wss_0$
obtained by a sequence of blow--ups of $S_0$ as before. We will call them ``main fibres''
since, in Section 3, we will prove that they are
 the only ones giving, by contractions,  all the possible singularities on a surface $S$ 
ruled by conics.

\medskip
\noindent
{\bf Notation.} Consider the above sequence of blow--ups.
If $D$ is a bisecant divisor on $S_0$, then we denote by $\wdd$ the
strict transform of $D$ on $\wss_0$ as well as on each surface $S_i$ defined in $(2)$.
While the strict transform of a component $e$ of a fibre
will be denoted by $\tilde e$ and also (for simplicity) by $e$ at each step of the sequence.

\medskip
\noindent
{\bf Remark 1.4.}
Let $\wdd \subset \wss_0$ be a bisecant divisor without base locus. Clearly,
$\Phi_{\wdd}$ is an isomorphism out of a finite number of fibres and assume, for
simplicity, that this number is one and this fibre is $f_n$. Keeping the notation
introduced in $(3)$, we then have:
$$
\harrowlength=40pt
\Phi_{\wdd}: \quad \wss_0 \setminus f_n \quad \mapright^{\cong} \quad S \setminus F.
$$
Therefore $\Phi_{\wdd}$ contracts $f_n$ to $F$. 
It is also clear that $\Phi_{\wdd}$ maps the divisor $\wdd$ to the hyperplane divisor
$H_{\wdd}$ of $S$, which is a surface ruled by conics. Hence $H_{\wdd}$ meets each fibre
of $S$ in two points; in particular, $H_{\wdd} \cdot F = 2 = \wdd \cdot f_n$.

\medskip
\noindent
{\bf Definition 1.5.} If $\wdd \subset \wss_0$ is as before, then we call
the ``$\wdd$--{\it degree}'' of a component $e$ of $f_n$ 
the integer 
$$
\deg_{\wdd}(e) = \wdd \cdot e.
$$

\goodbreak

\medskip
The arguments in 1.4 can be easily generalized to the case of non--empty base locus $B$
of
$|\wdd|$, giving immediately the following:

\medskip
\proclaim 
Proposition 1.6.  Let $\wss_0$, $\wdd$, $B$,  $f_n$ be as before and $e$ be
an irreducible component  of $f_n$. Then:
\item{$i)$}  $\wdd \cdot f_n =2$;
\item{$ii)$} if $e \not \subset B$ and $\deg_{\wdd}(e) = 0$ then $\Phi_{\wdd}$ contracts
$e$ to a point of the fibre $F$ of 
$S$;
\item{$iii)$} if $\deg_{\wdd}(e) < 0$  then $e\subset B$.
\cvd
\par

\medskip
\noindent
{\bf Remark 1.7.} Throughout this section, devoted to the construction of
``main'' fibres, we will assume that the centre of each blow--up 
$
\sigma_{P_i}: \; S_i \longrightarrow S_{i-1}
$
is a point $P_i \in S_{i-1}$ 
\item{-} either belonging to $\wdd \subset  S_{i-1}$ or infinitely near to it, i.e. $P_i$
belongs to the total transform of $D$ in $S_{i-1}$;
\item{-} belonging to a component of positive $\wdd$--degree of $f_{i-1} \subset
S_{i-1}$.
\par
\noindent
It is easy to see that the first condition is a consequence of the second one.

\medskip
\noindent
{\bf Remark 1.8.} Let $\sigma_{P_i}: \; S_i \longrightarrow S_{i-1}$ be as before and
let $g_{i-1}$ be a component of $f_{i-1} \subset S_{i-1}$ containing $P_i$. Denote as usual
by $\tilde g_{i-1}$ the strict transform of $g_{i-1}$ in $S_i$ and by $e$
the exceptional divisor of the blow--up. Then $\deg_{\wdd}(\tilde g_{i-1})=
\deg_{\wdd}(g_{i-1}) -1$ and $\deg_{\wdd}(e)=1$. In other words
$$
\deg_{\wdd}(\tilde g_{i-1})+ \deg_{\wdd}(e)=
\deg_{\wdd}(g_{i-1}) .
$$

\medskip
\noindent
{\bf Notation.} In the sequel we will draw the pictures of the fibres using  the
following agreement. Let $f_i \subset S_i$ be as before and let $e$ be one of its
irreducible components.
\item {-}
  The self--intersection of $e$ is ``represented''
as its degree i.e. if  $e^2 = -1, -2, -3 \dots$ then $e$ will be drawn as a line, a
conic, a cubic, etc. respectively. 
\item {-} The $\wdd$--degree of  $e$  is represented  by a continous line if
$\deg_{\wdd}(e) = 1$, a dashed line if $\deg_{\wdd}(e) = 0$ and a dotted line if
$\deg_{\wdd}(e) < 0$.  
\item {-}  The fibre $f_0$ of $S_0$ will be represented
as a continous smooth conic, even if $f_0^2 = 0$ and
$\deg_{D}(f_0) = 2$. 
\par

\medskip
\noindent
{\bf Remark 1.9.} Let us describe the fibres of level 0, 1, 2 arising from  blow--ups satisying the conditions in 1.7.

\smallskip
\noindent
{\bf Level 0.} By definition there is only one fibre of level 0 on $\wss_0 = S_0$, which
is $f_0$ itself.

\smallskip
\noindent
{\bf Level 1.} In this case also, there is only one fibre $f_1$ of $S_1 = \wss_0$
obtained from
$f_0$ by blowing--up
$S_0$ at a point $P_1 \in f_0$. Denoting by $e_1$ the exceptional divisor, it is clear
that $f_1 = f_0 + e_1$ and
$e_1^2 = -1$. Since $f_1^2 = 0$, it is immediate to see that also the other component
$f_0$ has self--intersection $-1$. Moreover, by 1.8,   
$\deg_{\wdd}(f_0) = \deg_{\wdd}(e_1) =1$.

  \medskip
\noindent
{\bf Level 2.} In order to get a fibre of level $2$, we can blow--up $f_1$ in a point
$P_2$ which is either a smooth point (i.e. belonging to exactly one of the two
components) or the singular point of $f_1$. These two cases are deeply different, so we
denote the corresponding fibres in a different way: $f_2(A)$ and
$f_2(D)$, respectively. Finally note that, in the case $(A)$, we can assume that $P_2$
belongs to one specific component, since the other construction can be recovered from
this one by an elementary transformation.

\medskip
\noindent
{\bf Case $(A)$.} Assume that $P_2 \in e_1$ and consider $\sigma_{P_2}: S_2 \longrightarrow S_1$. The fibre
$f_2$ of $S_2$ consists of three components: $f_2 = f_0 + e_1 + e_2$, where $e_2$ is
the exceptional divisor of $P_2$. From 1.8,
it is clear that $\deg_{\wdd}(e_1) = 0$, while 
$\deg_{\wdd}(f_0) = \deg_{\wdd}(e_2) = 1$. On the other hand, the self--intersection of
$e_2$ is $-1$, while the one of $e_1$ drops by one; briefly: $f_0^2 = e_2^2 = -1$ and
$e_1^2 = -2$. We can then draw the picture of
$f_2(A)$ using  the agreement introduced before.

\medskip    
   \centerline{
 \epsfxsize=7cm 
  \epsfbox{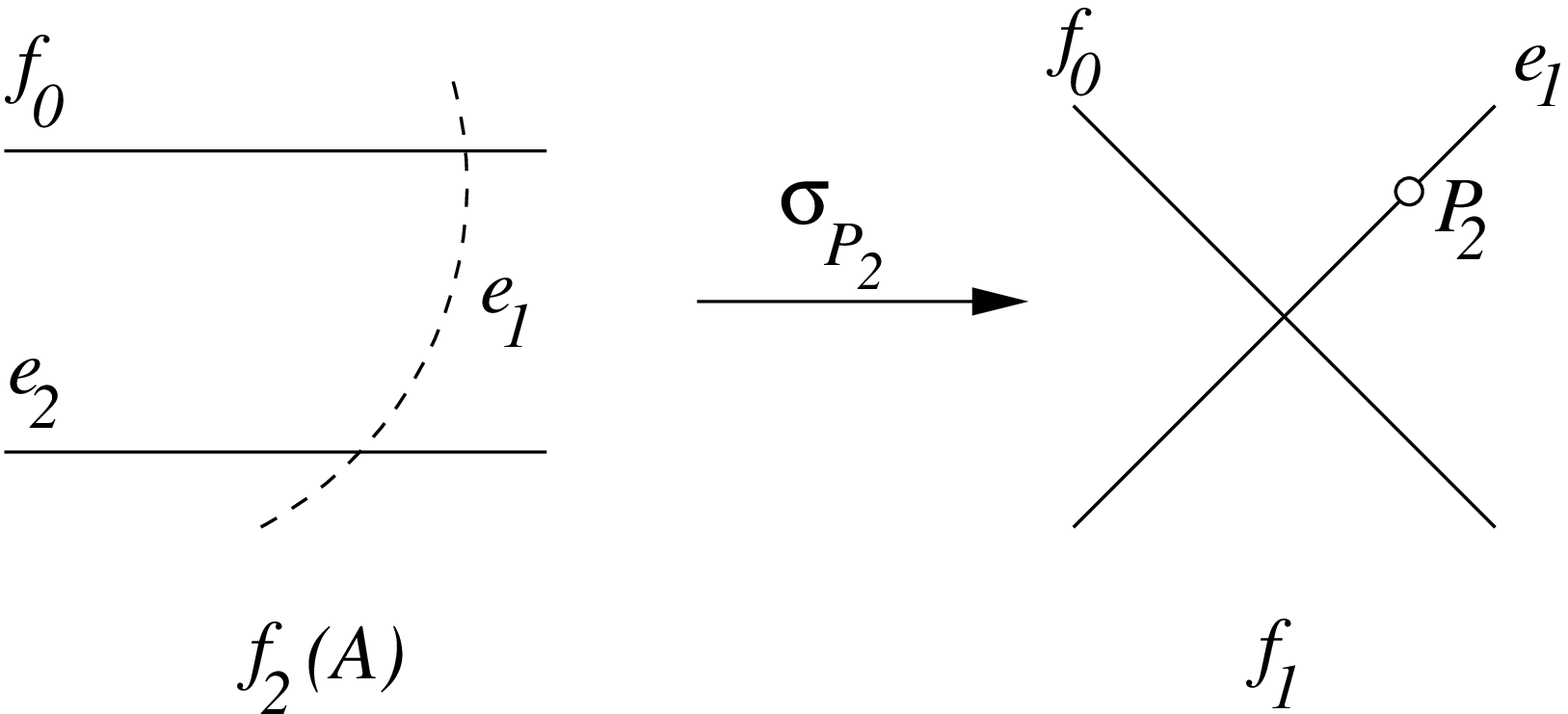}
 }
 \centerline {Figure 1}
 \medskip

\goodbreak

\noindent
{\bf Case $(D)$.} Assume that $P_2 = f_0 \cdot e_1$ and consider $\sigma_{P_2}: S_2
\longrightarrow S_1$. Now the exceptional divisor consists of a component of
multiplicity $2$ since $P_2$ is a double point of
$f_1$. Therefore the fibre $f_2$ of $S_2$ has the form $f_2 = f_0 + e_1 + 2e_2$. 
Since $\wdd$ is still bisecant on the fibre $f_2$, from 1.8 we have:
$\deg_{\wdd}(f_0) =
\deg_{\wdd}(e_1) = 0$ and $\deg_{\wdd}(2e_2) = 2$. On the other hand 
$f_0^2 = e_1^2 = -2$ and $e_2^2 = -1$, as  the following picture  shows.
 \medskip    
   \centerline{
 \epsfxsize=5cm 
\epsfbox{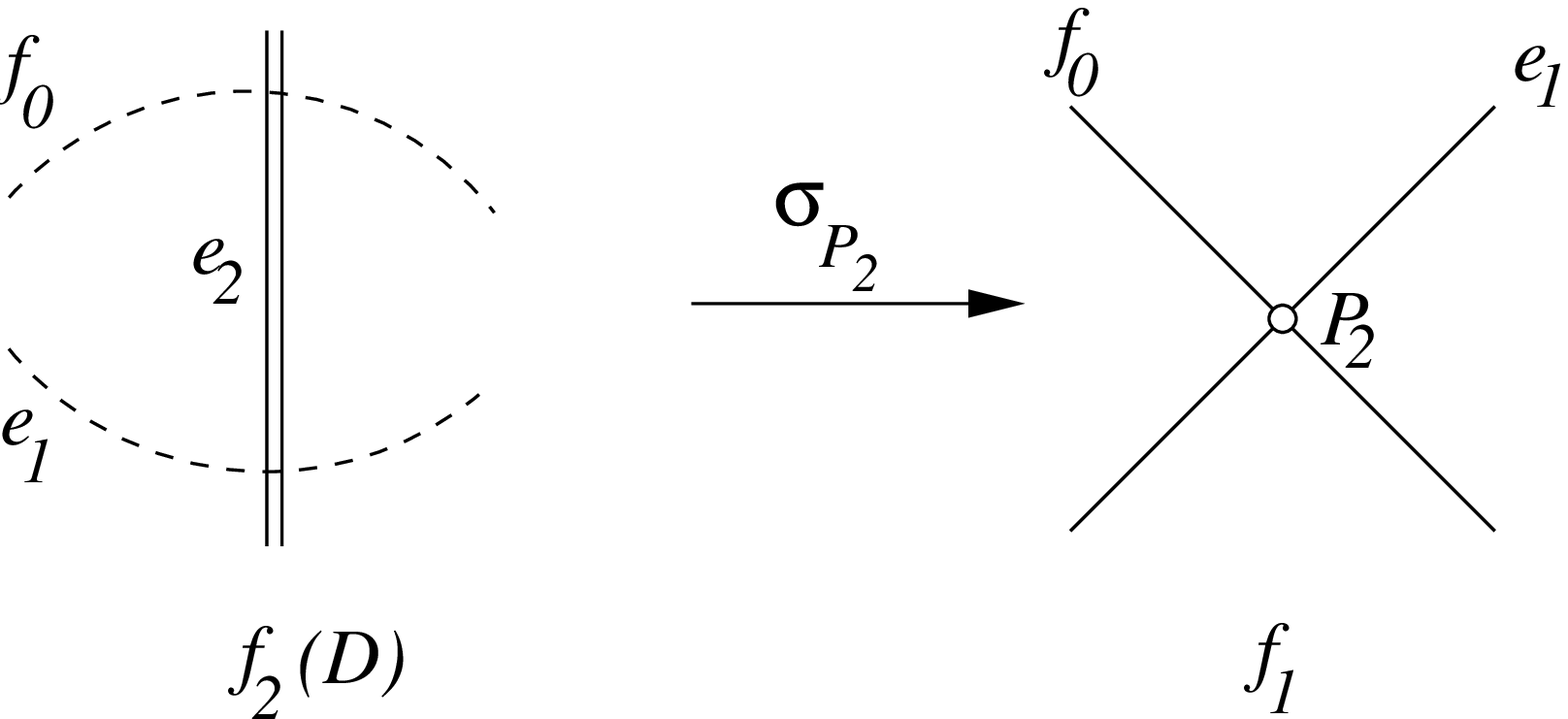}
 }
 \centerline {Figure 2}
 \smallskip

 Note that, in all the previous cases, there are no components with
negative $\wdd$--degree, hence 1.7 is automatically fulfilled. Moreover $|\wdd|$ is
base--point--free both on
$S_1$ and on
$S_2$.

\medskip
\noindent
{\bf Definition 1.10.} We define a  {\it main fibre} of level $n$ on $f_0$ recursively: 
\item {-} $f_1$ is the {\it main fibre of level $1$ on $f_0$};
\item {-}  if $n \ge 1$, then $f_n$  is a {\it main fibre of level $n$ on $f_0$} if it is
a fibre of $S_n = Bl_{P_n}(S_{n-1})$, where $P_n \in f_{n-1}$ and $f_{n-1}$ is a main
fibre of level $n-1$ on $f_0$. 
\item {-} if $n \ge 3$ then $P_n$ does not belong to  components of  $\wdd$--degree
$\le 0$.

\medskip
\noindent
{\bf Level $n$.} Now it is clear how to iterate the above construction. The ``main fibres'' of level $n$,
for $n \ge 3$, either of type $(A)$ and of type $(D)$ have the shapes described in the following picture:

 \medskip    
   \centerline{
 \epsfxsize=11cm 
\epsfbox{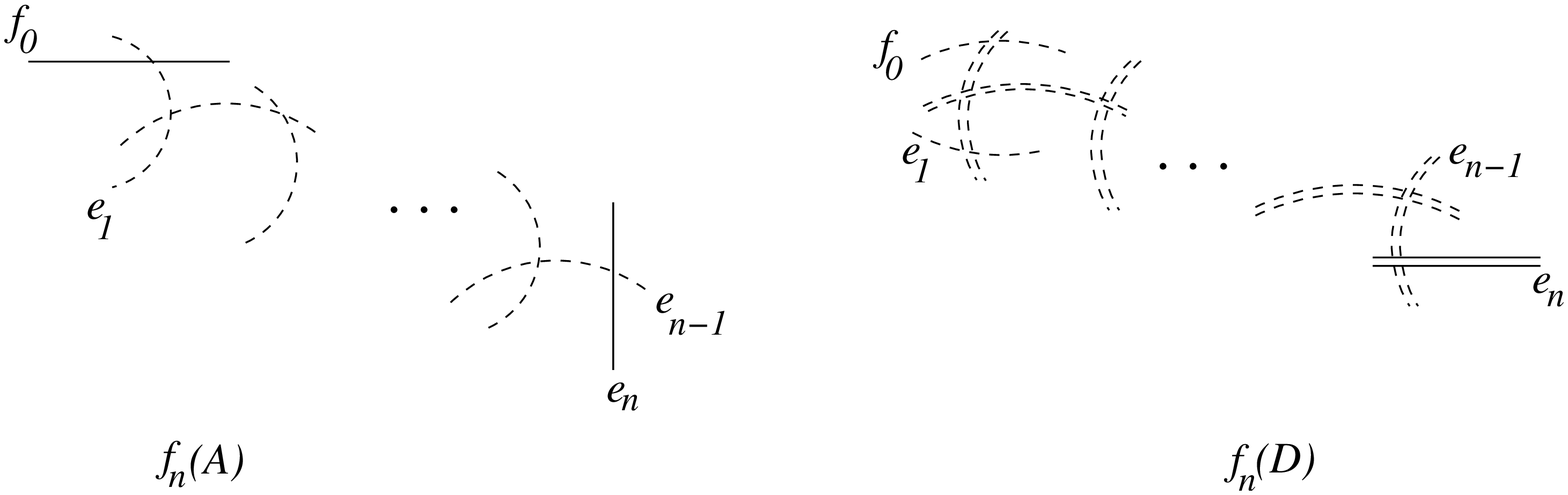}
 }
 \centerline {Figure 3}
 \smallskip

Let us describe the singularities of a surface $S$ ruled by conics as in $(3)$, 
where $\wss_0 = Bl_{P_1, \dots, P_n} (S_0)$,  the point $P_1$ belongs to a fibre $f_0$ of $S_0$ and  $P_2, \dots, P_n$ are such that the
corresponding fibre
 of $\wss_0$ is a main fibre $f_n$ of level $n$. Denote by $F_n$ the corresponding fibre on $S$, i.e.
$F_n := \Phi_{\wdd -B}(f_n)$.

\medskip
\noindent 
{\bf Remark 1.11.} Since $f_n$ is a main fibre, it is clear that all its components  have
$\wdd$--degree $0$, but one or two, having $\wdd$--degree $2$ or $1$, respectively (see
Figure 3). Moreover, it is easy to see that the above construction of main fibres leads
to a linear system $| \wdd |$ which is base--point--free.
 Finally, as noted in 1.6, the morphism
$\Phi_{\wdd}$ contracts exactly the components of $\wdd$--degree zero.

\medskip
\noindent 
{\bf Level 1.} The fibre  is $f_1 = f_0 + e_1$ and both
have $\wdd$--degree one. So  $\Phi_{\wdd}$ is an
isomorphism and  $F_1 \cong f_1$. 

\smallskip
\noindent 
{\bf Level 2.} 
{\bf Case $(A)$.} There is only one component of $\wdd$--degree zero,
namely $e_1$. Hence  $\Phi_{\wdd} = con (e_1)$. Since $e_1^2 = -2$, then the fibre $F_2(A)$ consists of
$f_0 + e_2$ and the point $f_0 \cap e_2$ is an ordinary double point of $S$. \hb
{\bf Case $(D)$.} Here two components have $\wdd$--degree zero, hence $\Phi_{\wdd} = con (f_0,e_1)$, and 
$f_0^2 = e_1^2 = -2$. Therefore the fibre $F_2(D)$ is $2e_2$ and its points $Q_0:= \Phi_{\wdd}(f_0)$, 
$Q_1:= \Phi_{\wdd}(e_1)$ are ordinary double points of $S$.

 \medskip    
   \centerline{
 \epsfxsize=11cm 
\epsfbox{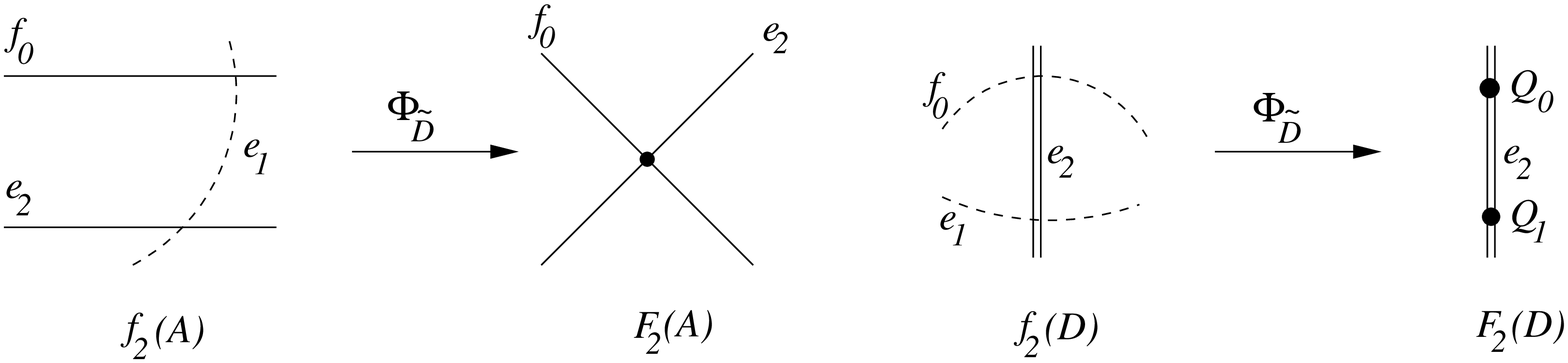}
 }
 \centerline {Figure 4}

\goodbreak

\medskip
\noindent 
{\bf Level 3.} \hb
{\bf Case $(A)$.} As before, since there are 2 components of $\wdd$--degree zero,  $\Phi_{\wdd} = con
(e_1, e_2)$ and $F_3(A) = f_0 + e_3$, where $Q:=f_0 \cap e_3$ is a rational double point of $S$. \hb
{\bf Case $(D)$.} In this case there are 3 components of $\wdd$--degree zero: $f_0, e_1, 2e_2$, so
$\Phi_{\wdd} = con (f_0, e_1, e_2)$ and $F_3(D) = 2e_3$, where $Q:=\Phi_{\wdd}(f_0+e_1+ 2e_2)$ is a
rational double point of $S$.

 \medskip    
   \centerline{
 \epsfxsize=13cm 
 \epsfbox{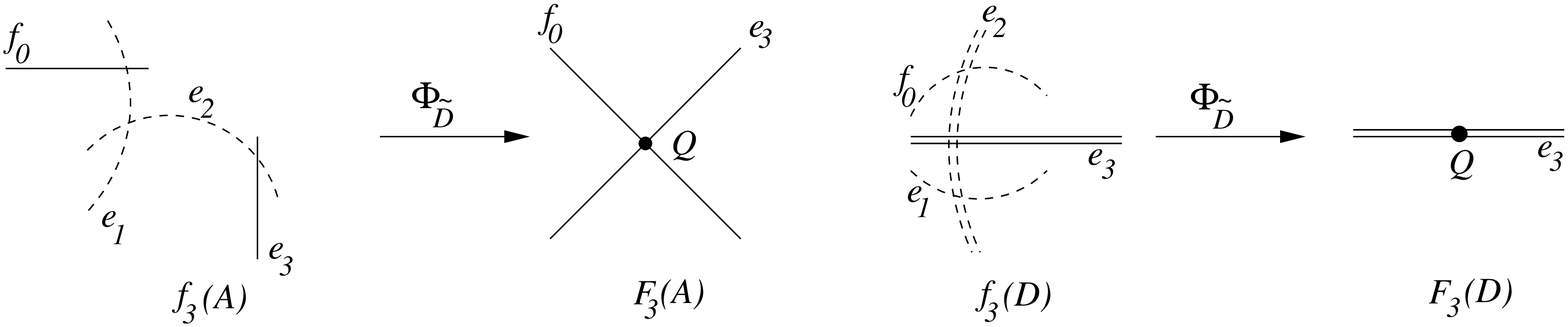}
 }
 \centerline {Figure 5}
 \medskip
 
\noindent 
{\bf Level $n$.} With the aid of the first cases and Figure 3,  the shapes of
the fibres $F_n$ (where $n \ge 3$) of $S$ and their singularities are clear. Namely the morphism
$\Phi_{\wdd}$ is exactly the blow--up of $S$ in a singular point and the union of the
components with $\wdd$--degree zero is the exceptional divisor of the blow--up.

\medskip
Let us briefly recall the notion  of  rational double point of a projective
surface.

\medskip
\noindent
{\bf Definition 1.12.} A {\it rational double point} $Q$ of a surface $S$ is  a double point such that
the exceptional set $E = E_1 + \cdots + E_n$  on the surface $\wss = Bl_Q(S)$ fulfills the
following requirements:
\item{$i)$} $E_1 \cup \cdots \cup E_n$ is connected; 
\item{$ii)$} $E_i \cong \Bbb P^1$ for all $i=1, \dots, n$; 
\item{$iii)$} $E_i^2 = -2$ for all $i=1, \dots, n$; 
\item{$iv)$} the lattice $\langle E_1, \dots, E_n \rangle$ generated by the irreducible components of $E$ is
negative definite.
\par

\medskip
\noindent
It is  well--known that  a rational double point of a surface is of one of the following types:
$A_n, D_n, E_6, E_7, E_8$ (see [1], 3.32). 
The dual graphs of the rational double points of type $A_n$ and $D_n$ are the following, where the $n$
components are represented by the $n$ vertices of the graph and two meeting components are connected by a
segment (see for instance [5], Ch.3):
\smallskip
$$
\circ \; \raise 2pt \hbox to 20 pt{\hrulefill} \,
\circ \, \raise 2pt \hbox to 20 pt{\hrulefill} \,
\circ \, \raise 2pt \hbox to 20 pt{\hrulefill} \;
\circ \; \cdots \cdots\; 
\circ \; \raise 2pt \hbox to 20 pt{\hrulefill} \,  
 \, \circ \; \raise 2pt \hbox to 20 pt{\hrulefill} \,  \circ
\eqno{A_n, \; n \ge 1}
$$

$$
\eqalign{
\circ \; \raise 2pt \hbox to 20 pt{\hrulefill} \,
\circ \, \raise 2pt \hbox to 20 pt{\hrulefill} \,
\circ \, \raise 2pt \hbox to 20 pt{\hrulefill} \;
\circ \; \cdots \cdots\; 
\circ \; \raise 2pt \hbox to 20 pt{\hrulefill} \,  
& \, \circ \; \raise 2pt \hbox to 20 pt{\hrulefill} \,  \circ \cr
& \, \, \  {\vrule height20pt} \cr
& \; \,  \raise 10pt \hbox{$\circ$}\cr
}
\eqno{D_n, \; n \ge 4} 
$$
Note that an $(A_1)$--singularity is an ordinary double point.

\medskip
We can summarize the above study of the singularities of $S$ arising from main fibres:

\proclaim
Proposition 1.13. Let $S_0= \Phi_D(\Bbb P ({\cal F}))$ be a surface geometrically ruled by
conics, 
$\harrowlength=25pt \wss_0 \mapright^\sigma S_0$ be a blow--up which is an isomorphism out of the fibre $f_n \subset
\wss_0$ and assume that $f_n$ is a main fibre of level $n$ over $f_0 \subset S_0$. Then $\wdd =
\sigma^{*}(D)$ is base--point--free and the morphism $\Phi_{\wdd}: \wss_0 \longrightarrow S$ is an
isomorphism on $\wss_0 \setminus f_n$. Moreover the fibre $F_n := \Phi_{\wdd}(f_n)$ of $S$ is of one of the
following types: \hb
$\bullet \; n=1:$ $F_1$ is the union of two distinct lines and $S$ is smooth (in this case $\Phi_{\wdd}$ is an
isomorphism everywhere);
\smallskip 
\noindent
$\bullet \; n=2:$ $F_2(A)$ is the union of two distinct lines, whose common point is an ordinary double
point of $S$; 
\item{}
 $\qquad F_2(D)$ is the union of two coincident lines, containing exactly two ordinary double
points of $S$;
\smallskip 
\noindent
$\bullet \; n \ge 3;$  $F_n(A)$ is the union of two distinct lines, meeting in a rational
double point of type $(A_{n-1})$;
\item{} $\qquad F_n(D)$ is  the union of two coincident lines, containing exactly
one rational double point of $S$; in particular, this point is of type $(A_3)$, if $n=3$, and of type
$(D_n)$, if $n \ge 4$.
\cvd
\par

\bigskip
\noindent
{\bf 2. Singularities arising from elementary transformations of  any  type} 

\medskip
\noindent 
{\bf Remark 2.1.} All the fibres $F_n$ of $S$ arising from main fibres $f_n$ of $\wss_0
= S_n$ (described in 1.13) are embedded fibres of level $n$. Namely, let $F_n$ be as
before and assume that $f_m \subset S_m$ is a fibre of level
$m$ which gives rise to $F_n$. Since $S_m$ is smooth at each point of $f_m$ (from 1.1),
then $f_m$ has to be obtained from $F_n$ by at least $n$ blow--ups. Therefore $m \ge n$,
so (from definition 1.3) the level of the embedded fibre $F_n$ is exactly $n$.

\medskip
The purpose of this section is to show that a fibre $F \subset S$ arising from a {\it
not main} fibre can be obtained  from a suitable {\it main} fibre (possibly of different
level).

\medskip
As usual, let us begin by describing the first case.

\medskip
\noindent 
{\bf Example 2.2.} Consider the fibre $f_2(A)$. In order to obtain a {\it not main} fibre
$f_3$ we have to blow--up $S_2$ in a point belonging to a component of $f_2(A)$ having $D$--degree $\le
0$ (see 1.5), hence either in  the vertex
$P = e_1
\cap e_2$ (or equivalently
$f_0
\cap e_2$) or in a point $Q \in e_1$ which is not a vertex, as the following picture illustrates.

\medskip    
   \centerline{
 \epsfxsize=8cm 
\epsfbox{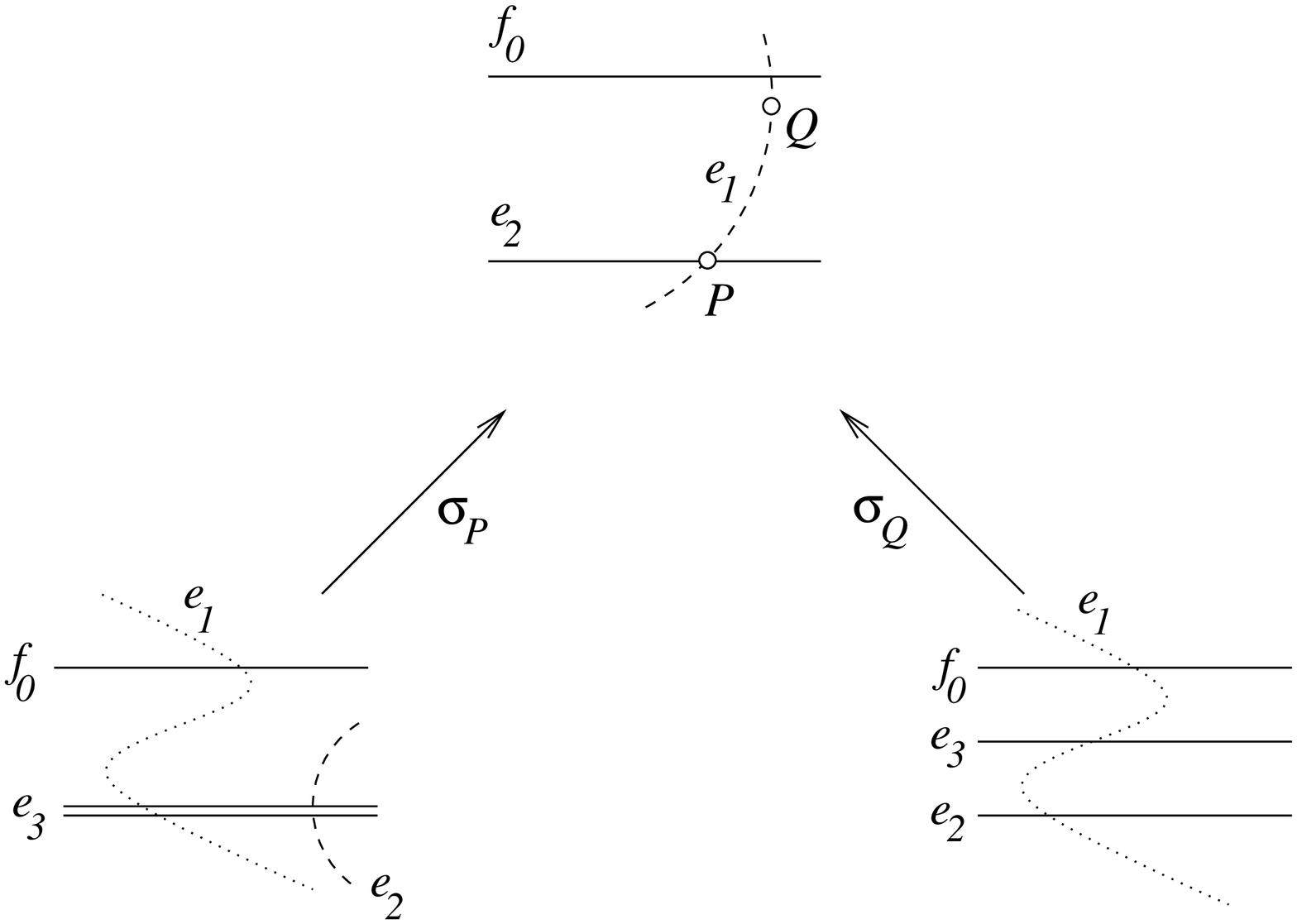}
 }
 \centerline {Figure 6}
 \medskip

In the first case, the fibre is $f_3:= f_0 +e_1 + e_2 + 2e_3$. Note that the intersections are described in
the above picture and the following relations hold:
$$
\displaylines{
f_0^2 = -1, \quad e_1^2 = -3, \quad
e_2^2 = -2, \quad e_3^2 = -1 \cr
\wdd \cdot f_0 = 1, \quad \wdd \cdot e_1 = -1, \quad
\wdd \cdot e_2 = 0, \quad \wdd \cdot e_3 = 1. \cr}
$$
Therefore $e_1 \subset B$. Then we have to compute the $(\wdd - e_1)$--degree of the
components:
$$
(\wdd - e_1) \cdot f_0 = 0, \quad (\wdd - e_1) \cdot e_1 = 2, \quad
(\wdd - e_1) \cdot e_2 = 0, \quad (\wdd - e_1) \cdot e_3 = 0.
$$
This proves that $\Phi_{\wdd - e_1}$ contracts all the components of $f_3$ but $e_1$,
hence $\Phi_{\wdd - e_1}(e_1)$ is a smooth conic: the  main fibre
of level $0$. \hb 
In the second case, the fibre is $f_3:= f_0 +e_1 + e_2 + e_3$. It can
be easily shown that:
$$
\displaylines{
f_0^2 = -1, \quad e_1^2 = -3, \quad
e_2^2 = -1, \quad e_3^2 = -1 \cr
\wdd \cdot f_0 = 1, \quad \wdd \cdot e_1 = -1, \quad
\wdd \cdot e_2 = 1, \quad \wdd \cdot e_3 = 1. \cr}
$$
As before, $e_1 \subset B$. So we have to compute the $(\wdd - e_1)$--degree of the components:
$$
(\wdd - e_1) \cdot f_0 = 0, \quad (\wdd - e_1) \cdot e_1 = 2, \quad
(\wdd - e_1) \cdot e_2 = 0, \quad (\wdd - e_1) \cdot e_3 = 0.
$$
As in the previous case, $\Phi_{\wdd - e_1}$ contracts $f_0,e_2,e_3$, while $\Phi_{\wdd - e_1}(e_1)$ is a
smooth conic, so isomorphic to $f_0$.

\medskip
The case $f_2(D)$ is analogous. Both cases  show that, if $f_3$ is any possible fibre of level $3$, but not a main fibre, then $F_3 =
\Phi_{\wdd - B}(f_3)$ is   of type $f_0$. 
Concerning the general case, we want to show the following:

\medskip
\noindent
{\bf Claim:}
Let $F \subset S$ be a fibre obtained as $\Phi_{\wdd - B}(f_n)$, where $f_n$ is a fibre
of level $n$ of $\wss_0 = S_n$. If $f_n$ is {\it not} a main fibre, then   $F$ is an
embedded fibre of level $m$, where $m <n$, i.e.  it can be obtained as $F= \Phi_{\wdd - B}(f_m)$, where
$f_m$ is a main fibre of level $m$ on $S_m$ for a suitable $m < n$.

\medskip
\noindent
{\bf Remark 2.3.} It i is enough to show the claim when $f_n$ is
not a main fibre, but  comes from a main fibre by a blow--up in a ``not admissible''
single point, in the sense of 1.10. So we can assume that:
\item {-} $f_{n-1}$  is a main fibre of level $n-1$ on $S_{n-1}$;
\item {-} $P_n \in f_{n-1}$; 
\item {-} $P_n$  belongs to a component of $f_{n-1}$ having $\wdd$--degree zero.

  \medskip    
   \centerline{
 \epsfxsize=12cm 
\epsfbox{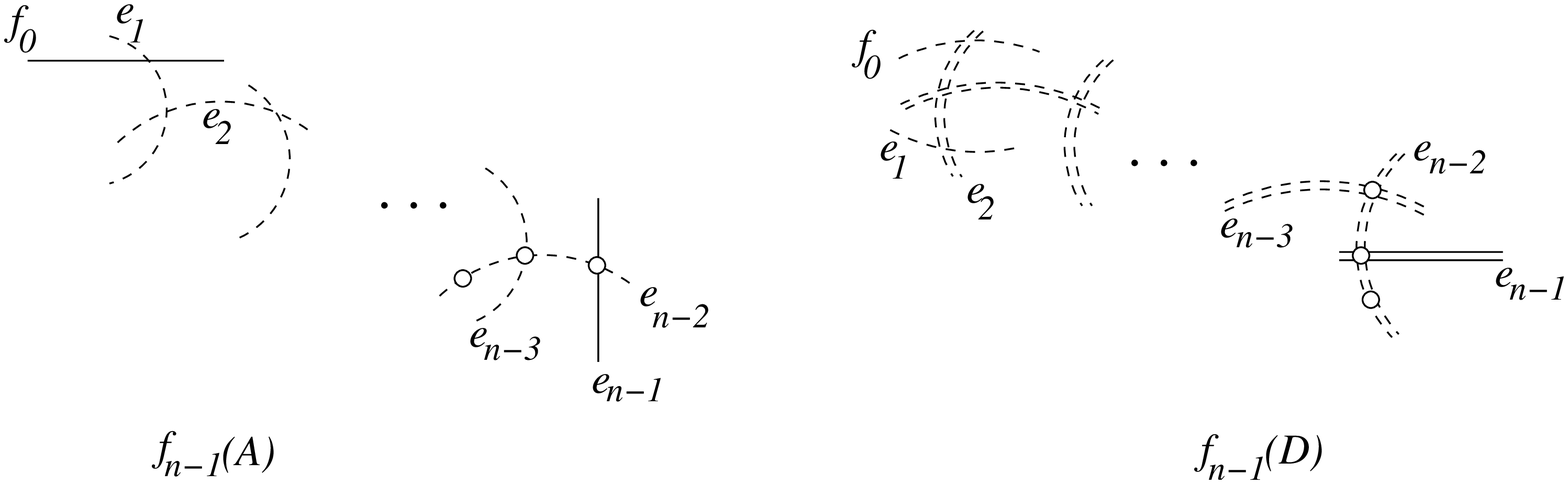}
 }
 \centerline {Figure 7}

 \medskip
\noindent
Looking at the picture above, it is clear that
the cases to be considered are the following: 
\smallskip
\item {-} $P_n \in f_{n-1}(A)$
\item\item {$(1)$} $P_n$ is  a vertex and belongs also to a component of self--intersection $-1$ (e.g. $P_n
=e_{n-1} \cap  e_{n-2}$); 
\item\item {$(2)$}  
$P_n$ is not a vertex (e.g. $P_n \in e_{n-2}$ and no other component);
\item\item {$(3)$}  $P_n$ is  a vertex belonging to two components having both self--intersection $-2$ (e.g.
$P_n =e_{n-2} \cap  e_{n-3}$).
\smallskip
\item {-} $P_n \in f_{n-1}(D)$
\item\item {$(4)$}  $P_n$ is  a vertex and belongs also to a component of self--intersection $-1$ (e.g. $P_n
=e_{n-1} \cap  e_{n-2}$);
\item\item {$(5)$}  $P_n$ is not a vertex  (e.g. $P_n \in e_{n-2}$ and no other component); 
\item\item {$(6)$}  $P_n$ is  a vertex belonging to two components having both self--intersection $-2$ (e.g.
$P_n =e_{n-2} \cap  e_{n-3}$).
\par

\medskip
In order to compute the singularities arising in the above cases, let us introduce some notation. \hb
Let $D$ be the bisecant divisor on the surface $S_{n-1}$, so $\deg_D(f_{n-1}) = 2$. 
Let $C_0 := f_0$,  $C_i:=e_i$, for $i=1, \dots, n-1$. 
Set also $I_{n-1}:=(C_i \cdot C_j)_{i,j=0, \dots, n-1}$ the intersection matrix of the
main fibre $f_{n-1}$ of level
$n-1$ and $D_{n-1}$ the vector of the $D$--degrees of the components, i.e. $D_{n-1}:= (\deg_D(C_0),
\deg_D(C_1), \dots, \deg_D(C_{n-1}))$. \hb
Finally,  let us introduce the vector consisting
of the multiplicities of the components in the fibre $f_{n-1}$ (and of the other fibres arising from the
blow--ups). Set 
$$
\muu(f_{n-1}):= (\mu(C_0), \mu(C_1), \dots, \mu(C_{n-1})).
$$
Looking at Figure 7, it is clear that:
$\muu(f_{n-1}(A))= (1,1, \dots, 1))$ and  $\muu(f_{n-1}(D))= (1,1, 2, \dots, 2))$.
Let us begin from the data of the two fibres $f_{n-1}(A)$ and $f_{n-1}(D)$: 
$$
\underbrace{
\left | \matrix{
1 \cr
1 \cr
1 \cr
1 \cr
\cr
\vdots \cr
\cr
1 \cr
1 \cr
1 \cr
}
\right |
}_{\muu(f_{n-1}(A))}
\quad
\underbrace{\pmatrix{
-1 & 1 & 0 & 0 & \cdots & 0 & 0 & 0 \cr
1 & -2 & 1  & 0 & \cdots & 0 & 0 & 0 \cr
0 & 1 & -2  & 1 & \cdots & 0 & 0 & 0 \cr
0 & 0 & 1  & -2 & \cdots & 0 & 0 & 0 \cr
\cr
\vdots & &  &  &\ddots & & & \vdots \cr
\cr
0 & 0 & 0  & 0 & \cdots & -2 & 1 & 0 \cr
0 & 0 & 0  & 0 & \cdots & 1 & -2 & 1 \cr
0 & 0 & 0  & 0 & \cdots & 0 & 1 & -1 \cr
}
}_{I_{n-1}(A)}
\quad
\underbrace{
\left | \matrix{
1 \cr
0 \cr
0 \cr
0 \cr
\cr
\vdots \cr
\cr
0 \cr
0 \cr
1 \cr
}
\right |
}_{D_{n-1}(A)}
$$
and
$$
\underbrace{
\left | \matrix{
1 \cr
1 \cr
2 \cr
2 \cr
\cr
\vdots \cr
\cr
2 \cr
2 \cr
2 \cr
}
\right |
}_{\muu(f_{n-1}(D))}
\underbrace{\pmatrix{
-2 & 0 & 1 & 0 & \cdots & 0 & 0 & 0 \cr
0 & -2 & 1  & 0 & \cdots & 0 & 0 & 0 \cr
1 & 1 & -2  & 1 & \cdots & 0 & 0 & 0 \cr
0 & 0 & 1  & -2 & \cdots & 0 & 0 & 0 \cr
\cr
\vdots & &  &  &\ddots & & & \vdots \cr
\cr
0 & 0 & 0  & 0 & \cdots & -2 & 1 & 0 \cr
0 & 0 & 0  & 0 & \cdots & 1 & -2 & 1 \cr
0 & 0 & 0  & 0 & \cdots & 0 & 1 & -1 \cr
}
}_{I_{n-1}(D)}
\underbrace{
\left | \matrix{
0 \cr
0 \cr
0 \cr
0 \cr
\cr
\vdots \cr
\cr
0 \cr
0 \cr
1 \cr
}
\right |
}_{D_{n-1}(D)}
$$

\medskip
Note that, in both cases $(A)$ and $(D)$, the scalar product of $\muu(f_{n-1})$ and each row of
$I_{n-1}$ is zero, while $\muu(f_{n-1}) \cdot D_{n-1} = 2$. 

\medskip
Consider now the blow--up at $P_n$ and let  $\wdd$ be the strict transform of $D$ on
$S_n$. Let us denote again by $C_0, \dots, C_{n-1}, C_n$ the components of the fibre
$f_n$, where $C_n$ is the exceptional divisor of the blow--up $\sigma_{P_n}$. We denote 
by $\wii_n$ the intersection matrix of $f_n$ and by $\deg_{\wdd}$ the vector of the
$\wdd$--degrees of its components, i.e. 
$\deg_{\wdd}:= (\deg_{\wdd}(C_0),
\deg_{\wdd}(C_1), \dots, \deg_{\wdd}(C_{n}))$.
Finally set, for each $h$ and $k$:
$\Sigma_h^k = \sum_{i=h}^{k}C_{i}$.

\medskip
For sake of brevity, we examine explicitly only the first case.

\medskip
\noindent
{\bf Case (1): $f_{n-1}(A) \ni P_n = C_{n-1} \cap C_{n-2}$} \hb
 Since $\deg_{\wdd}(C_{n-2}) =
-1$, then $C_{n-2}$ is contained in the base locus of $\wdd$. So we compute the degrees
of the components of $f_n$   with respect to $\wdd - C_{n-2}$. It is immediate to see
that $\deg_{\wdd-C_{n-2}}(C_{n-3}) = -1$, so also $C_{n-3}$ is contained in the base
locus of $\wdd$. We then iterate this computation up to a degree vector whose components
are all non--negative. 
This is the intersection matrix $\wii_{n}$ and the list of the degree
vectors:
$$
\underbrace{
\left | \matrix{
1 \cr
1 \cr
1 \cr
1 \cr
\cr
\vdots \cr
\cr
1 \cr
1 \cr
1 \cr
1 \cr
1 \cr
2 \cr
}
\right |
}_{\muu}
\underbrace{\pmatrix{
-1 & 1 & 0 & 0 & \cdots & 0 & 0 & 0 & 0 & 0 & 0\cr
1 & -2 & 1  & 0 & \cdots & 0 & 0 & 0 & 0 & 0 & 0\cr
0 & 1 & -2  & 1 & \cdots & 0 & 0 & 0 & 0 & 0 & 0\cr
0 & 0 & 1  & -2 & \cdots & 0 & 0 & 0 & 0 & 0 & 0\cr
\cr
\vdots & &  &  &\ddots & & & & & & \vdots \cr
\cr
0 & 0 & 0  & 0 & \cdots & -2 & 1 & 0 & 0 & 0 & 0\cr
0 & 0 & 0  & 0 & \cdots & 1 & -2 & 1 & 0 & 0 & 0\cr
0 & 0 & 0  & 0 & \cdots & 0 & 1 & -2 & 1 & 0 & 0\cr
0 & 0 & 0  & 0 & \cdots & 0 & 0 & 1 & -3 & 0 & 1\cr
0 & 0 & 0  & 0 & \cdots & 0 & 0 & 0 & 0 & -2 & 1\cr
0 & 0 & 0  & 0 & \cdots & 0 & 0 & 0 & 1 & 1 & -1\cr
}
}_{\wii_{n}}
\;
\underbrace{
\left | \matrix{
1 \cr
0 \cr
0 \cr
0 \cr
\cr
\vdots \cr
\cr
0 \cr
0 \cr
0 \cr
-1 \cr
0 \cr
1 \cr
}
\right |
}_{\deg_{\wdd}}
\underbrace{
\left | \matrix{
1 \cr
0 \cr
0 \cr
0 \cr
\cr
\vdots \cr
\cr
0 \cr
0 \cr
-1 \cr
2 \cr
0 \cr
0 \cr
}
\right |
}_{(*_1)}
\underbrace{
\left | \matrix{
1 \cr
0 \cr
0 \cr
0 \cr
\cr
\vdots \cr
\cr
0 \cr
-1 \cr
1 \cr
1 \cr
0 \cr
0 \cr
}
\right |
}_{(*_2)}
\underbrace{
\left | \matrix{
0 \cr
1 \cr
0 \cr
0 \cr
\cr
\vdots \cr
\cr
-1 \cr
1 \cr
0 \cr
1 \cr
0 \cr
0 \cr
}
\right |
}_{(*_3)}
\; \cdots \;
\underbrace{
\left | \matrix{
1 \cr
-1 \cr
1 \cr
0 \cr
\cr
\vdots \cr
\cr
0 \cr
0 \cr
0 \cr
1 \cr
0 \cr
0 \cr
}
\right |
}_{(*_4)}
\underbrace{
\left | \matrix{
0 \cr
1 \cr
0 \cr
0 \cr
\cr
\vdots \cr
\cr
0 \cr
0 \cr
0 \cr
1 \cr
0 \cr
0 \cr
}
\right |
}_{(*_5)}
$$
where
$$
\matrix{
(*_1): \deg_{\wdd - C_{n-2}} & (*_2): \deg_{\wdd - \Sigma_{n-3}^{n-2}} & 
(*_3): \deg_{\wdd - \Sigma_{n-4}^{n-2}} & (*_4): \deg_{\wdd - \Sigma_{2}^{n-2}} & (*_5): \deg_{\wdd - \Sigma_{1}^{n-2}} \cr
}.
$$
\noindent
Hence set $B:=\Sigma_{1}^{n-2} = \sum_{i=1}^{n-2}C_{i}$. First note that $\Phi_{\wdd
-B}$ contracts $C_n$, hence the blow--up in $P_n$ is somewhat irrelevant and the image of
$\Phi_{\wdd -B}$  is a surface which can be recovered from a main fibre of lower level.
\hb
 In order to precisely understand the kind of the fibre we obtain, let us examine the
contracted locus. \hb 
As appears from the vector $\deg_{\wdd - \Sigma_{1}^{n-2}}$, it consists of three
connected components: $C_0$, $C_{n-1}
\cup C_n$ and $C_2 \cup \dots \cup C_{n-3}$. The first two components give rise to two
smooth points, respectively. Moreover, as the following picture shows, $\Phi_{\wdd - B}$
factorizes as 

\goodbreak

$$
\Phi_{\wdd - B} = con(C_0, C_2, C_3, \dots, C_{n-3},C_{n-1}, C_n) = 
con(C_2, C_3, \dots, C_{n-3}) \circ 
con(C_{n-1}) \circ
con(C_0, C_n).
$$

\medskip    
   \centerline{
 \epsfxsize=7.8cm 
\epsfbox{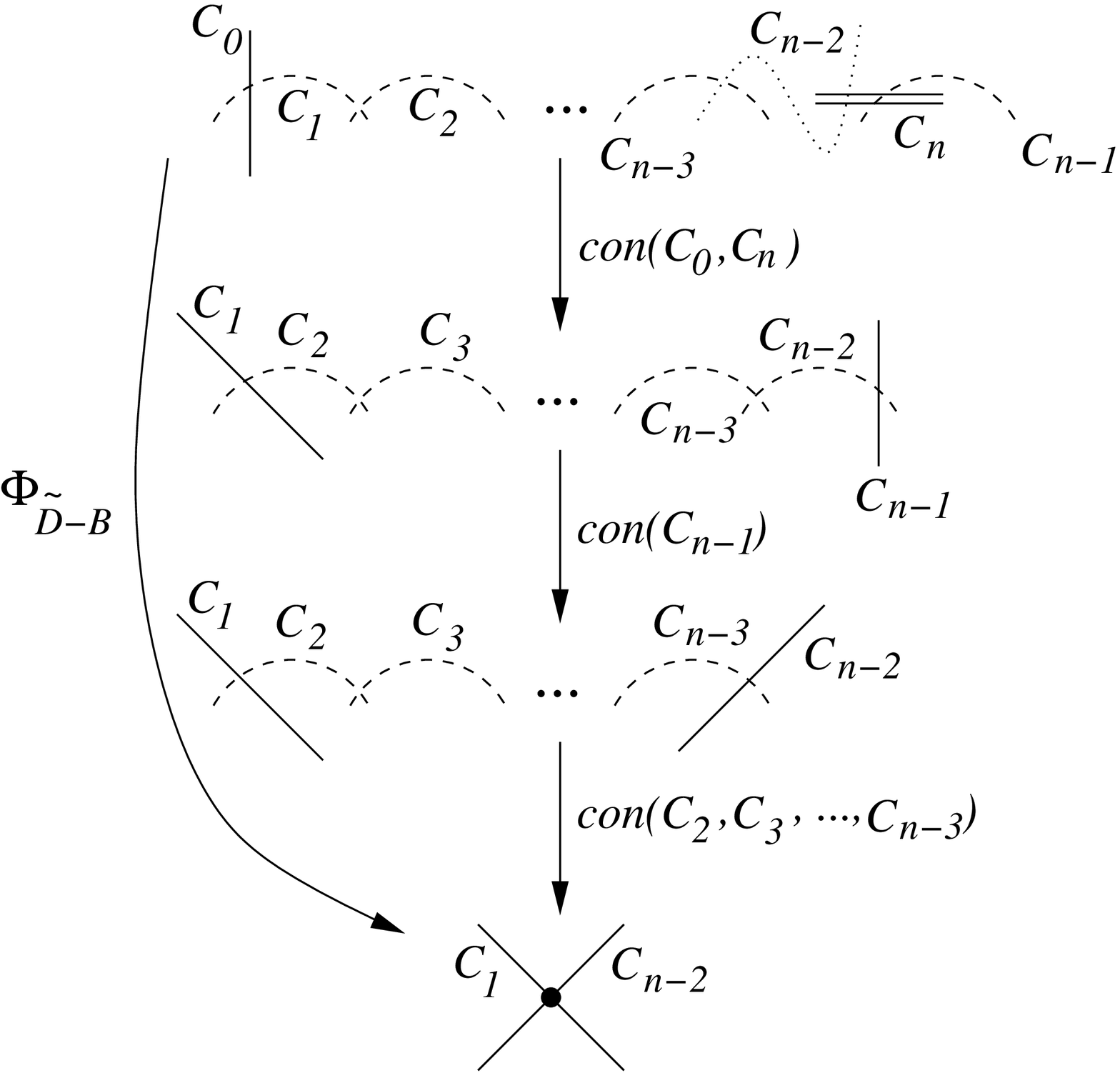}
 }
 \centerline {Figure 8}

\medskip
\noindent
On the other hand, the last component $C_2 \cup \dots \cup C_{n-3}$ of the contracted locus gives rise to a
rational double point of type
$(A_{n-4})$. Therefore
$F
\subset S =
\Phi_{\wdd -B}(\wss_0)$  is an embedded fibre of level $n-3$.

\bigskip
With similar arguments, the whole claim can be proved. Hence, using 1.13 and the claim, we have shown in this way the
following result:

\medskip
\proclaim 
Theorem 2.4. Let $S \subset \Bbb P^N$ be a projective surface ruled by conics over a smooth irreducible curve.
Then the degenerate fibres and the singular points of $S$ are as follows:
\item{$i)$} If a fibre  is the union of two distinct lines then $S$ is either smooth along this fibre or singular only at the common point of the lines; in this case the singularity is a rational double point of type $(A_n)$, $n \ge 1$.
\item{$i)$} A fibre which is the union of two coincident lines contains either exactly two ordinary double
points of $S$ or exactly one singular point of $S$; 
in this case the singularity is a rational double point either of type $(A_3)$ or of type $(D_n)$,  $n \ge 4$.
\cvd
\par

\medskip
\proclaim 
Corollary 2.5. Let $S \subset \Bbb P^N$ be a projective surface ruled by conics over a smooth irreducible curve.
Then the singular points of $S$ are rational double points of type $(A_n)$ or $(D_n)$.
\cvd

\medskip
\noindent
{\bf Remark 2.6.}
By 0.3 and 0.10, $(i) - ( ii)$, every projective surface $S$  ruled by conics can have at worst (finitely many) normal Gorenstein singularities. Assume, for sake of simplicity, that $S$ has exactly one singularity. Clearly, the exceptional divisor of the resolution of singularities $\wss_0$  of $S$ is contained in the degenerated fibre of the composition 
$$
\harrowlength=25pt \pi:= p \circ \sigma :  \; \wss_0 \mapright C
$$ 
where $\sigma : \wss_0 \rightarrow S$ is as in diagram $(1)$  and $p: \wss_0 = \Phi_D(\Bbb P({\cal F})) \longrightarrow C$ is the canonical projection.

Let us observe that a result of Badescu concerning nonrational ruled surfaces (see [1], Lemma 14.35) holds also in a slightly different situation: namely, if  $\pi:X \rightarrow B$  is a surjective morphism  from a smooth projective surface $X$ to a smooth rational curve $B$, whose general fibre is a smooth rational curve, then  $H^1({\cal O}_Z) = 0$ for every positive divisor $Z$ of support contained in a degenerated fibre of $\pi$ and such that the  intersection matrix of $Z$ is negative definite.

Hence, one can apply this result to the morphism $\pi :  \wss_0 \rightarrow C$, obtaining in particular that  $H^1({\cal O}_Z) = 0$ for every positive divisor $Z$ of support contained in the exceptional divisor of $\sigma : \wss_0 \rightarrow S$. 

Therefore, performing the same argument for each singular point of $S$,  one obtains via a criterion of M. Artin (see [1], Lemma 3.8 and Definition 3.17)  that every singularity of $S$ is rational.

Finally, (using [1], Corollary 4.19) any rational surface singularity which is Gorenstein is a rational double point. 

These general arguments show that  every projective surface ruled by conics can have only rational double points as singularities. However, the main result of this paper (beside the constructive procedure performed in sections 1 and 2) gives a  more precise description of the singularities: namely,  it tells us that every projective surface ruled by conics can have only rational double points of type ${\bf A}_n$  ($n \ge 1$), or of type ${\bf D}_n$  ($n \ge 4$), as singularities. In other words, the rational double points of type ${\bf E}_n$, with $n = 6, 7, 8$, cannot occur on such a surface.

\bigskip
\noindent
{\bf 3. Birational models of a surface ruled by conics}
\medskip
We want to describe the surfaces  geometrically ruled by conics which give rise
to a surface $S$ ruled by conics, i.e. we want to find the surfaces $S_0, S_0', \dots$ such that for each of them diagram 
$(1)$ holds, as follows:
$$
\commdiag{ 
S'_0& \mapleft^{\sigma'} &\wss_0 &  \mapright^\sigma & S_0 \cr
&&\mapdown \lft{\Phi_{\wdd -B}} \cr
&&S
}
$$
Since this is a local study, we can assume that $S$ has only one singular embedded fibre $F_n$ of level $n$. Clearly
 $F_n \subset S$  can be obtained in a unique way from a main fibre $f_n \subset \wss_0$
having the same level
$n$ (see 2.1). Therefore the initial problem can be reduced to the following: \hb
{\sl For each $n$, choose a main fibre $f_n$ of $\wss_0$. Describe all the surfaces
$S_0$ and the blow--ups $\sigma: \; \wss_0 \rightarrow S_0$ such that 
$f_n$ maps to a suitable fibre $f_0$ of $S$.}

\medskip
\noindent
{\bf Example 3.1.} The unique fibre of level $1$ is
$f_1 \subset \wss_0$, obtained from $f_0$ simply by blowing--up
$S_0$ at a point $P_1 \in f_0$. 
Clearly, $f_1 = f_0 \cup e_1$ can be obtained  in exactly one other way. If we  contract its component $f_0$, then we obtain a surface
$S'_0 = elm_{P_1}(S_0)$, whose fibre corresponding to
$f_1$ is a smooth conic $e_1$ and the component $f_0$ contracts to a (smooth) point,
say $F_0$, of $e_1$. \hb Finally, it is also clear that the above contraction is the
blow--up of $S'_0$ at $F_0$. The following picture describes this situation:

  \medskip    
   \centerline{
 \epsfxsize=12cm 
\epsfbox{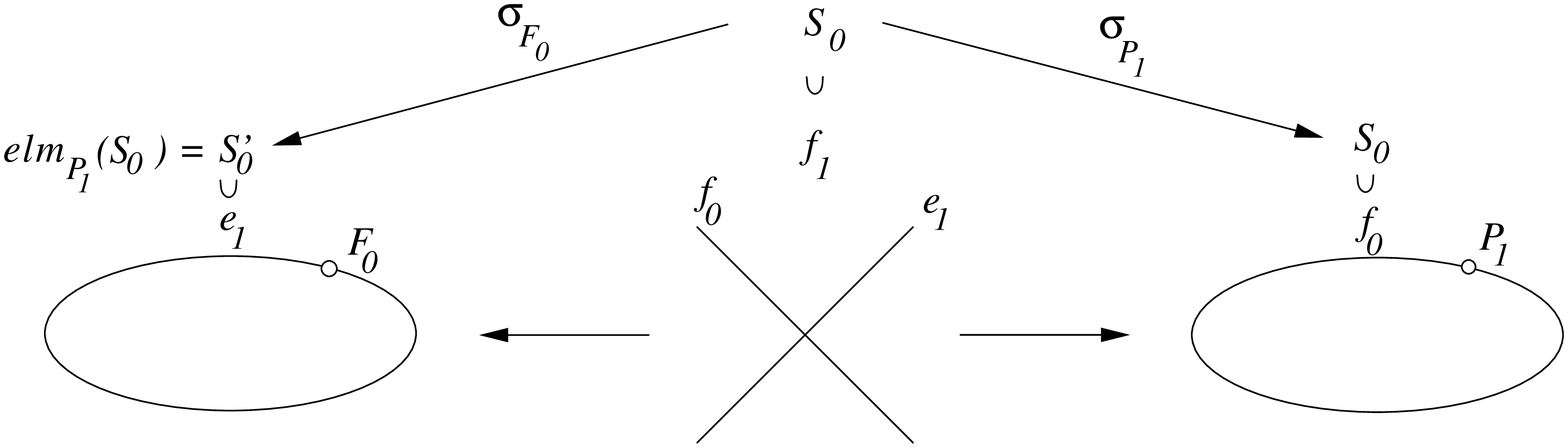}
 }
 \centerline {Figure 9}
 \medskip
 
\noindent
{\bf Example 3.2.} Consider now the main fibre $f_2(A) \subset \wss_0$, arising
from
$f_0$ via the blowing--up
$\sigma_{P_1P_2}$, where $P_1 \in f_0$ and $P_2$ is a point infinitely near to $P_1$ along a 
direction which is transversal to $f_0$. \hb
It is clear that $\sigma_{P_1P_2}$ contracts the components $e_1$ and $e_2$ of
$f_2(A)$. The following figure illustrates the other two surfaces, $S'_0$ and
$S''_0$ say, obtained by contracting $f_0$ and $e_2$ (in the middle) and  $f_0$ and $e_1$ (on the left hand side),  respectively.
 \medskip    
   \centerline{
 \epsfxsize=12cm 
\epsfbox{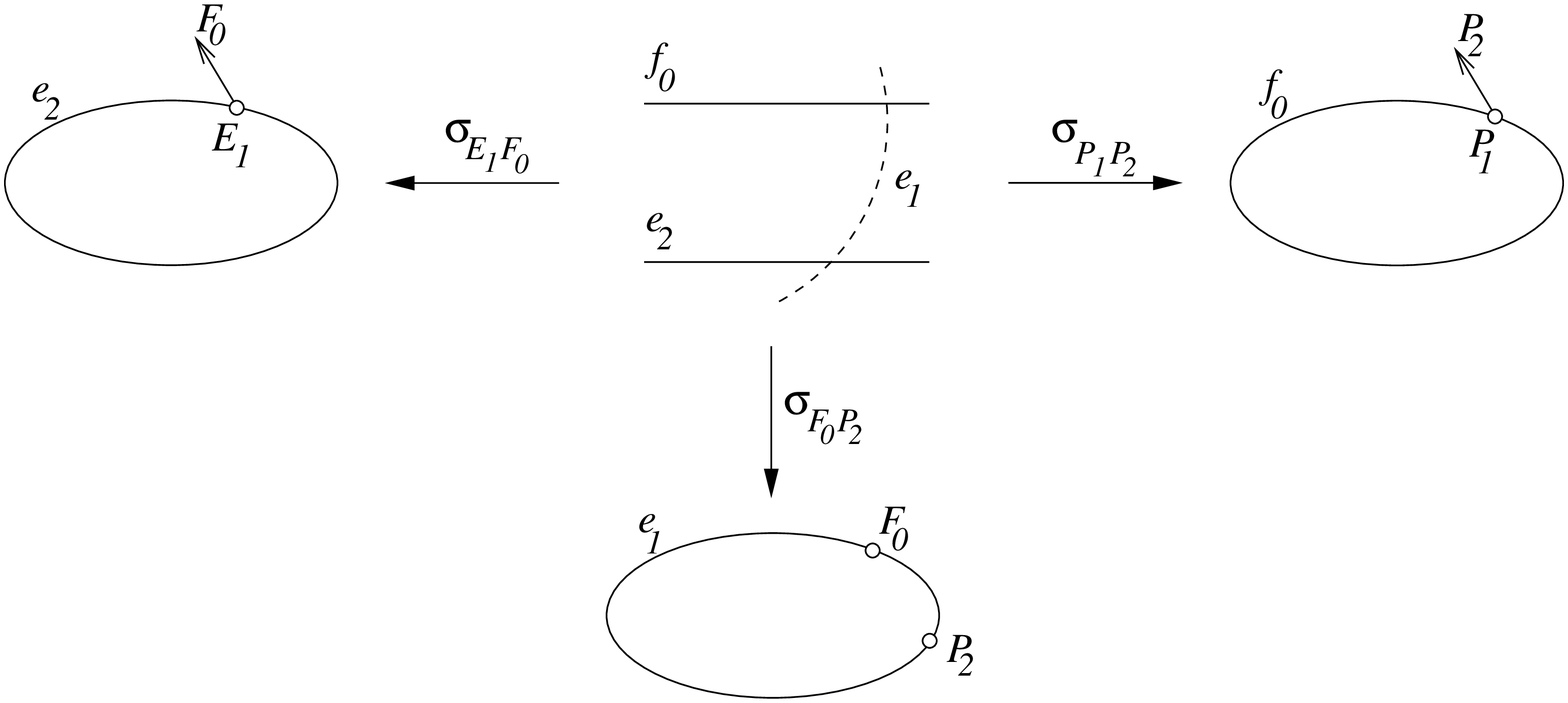}
 }
 \centerline {Figure 10}
 \medskip
 
\noindent
In the first case, we contract $f_0$ and $e_2$ to two distinct points, say $F_0$ and $E_2$,
of the fibre $e_1 \subset S'_0$. \hb
In the second case, we obtain a smooth conic $e_2 \subset S''_0$ and the contraction turns out
to be the blowing--up of $S''_0$ at $E_1 \in e_2$ and at $F_0$, a point infinitely near to
$E_1$. The exceptional divisor of $\sigma_{E_1}$ is $e_1$, moreover $F_0 \in e_1$ and its
exceptional divisor is $e_2$.   Hence $\sigma_{P_1P_2}= con(e_1,e_2)$, 
$\sigma_{F_0P_2}= con(f_0,e_2)$, 
$\sigma_{EF_0}= con(e_1,f_0)$. \hb
Finally note that the surfaces $S_0, S'_0, S''_0$ are
related by elementary transformations, as in the following picture:

 \medskip    
   \centerline{
 \epsfxsize=13cm 
\epsfbox{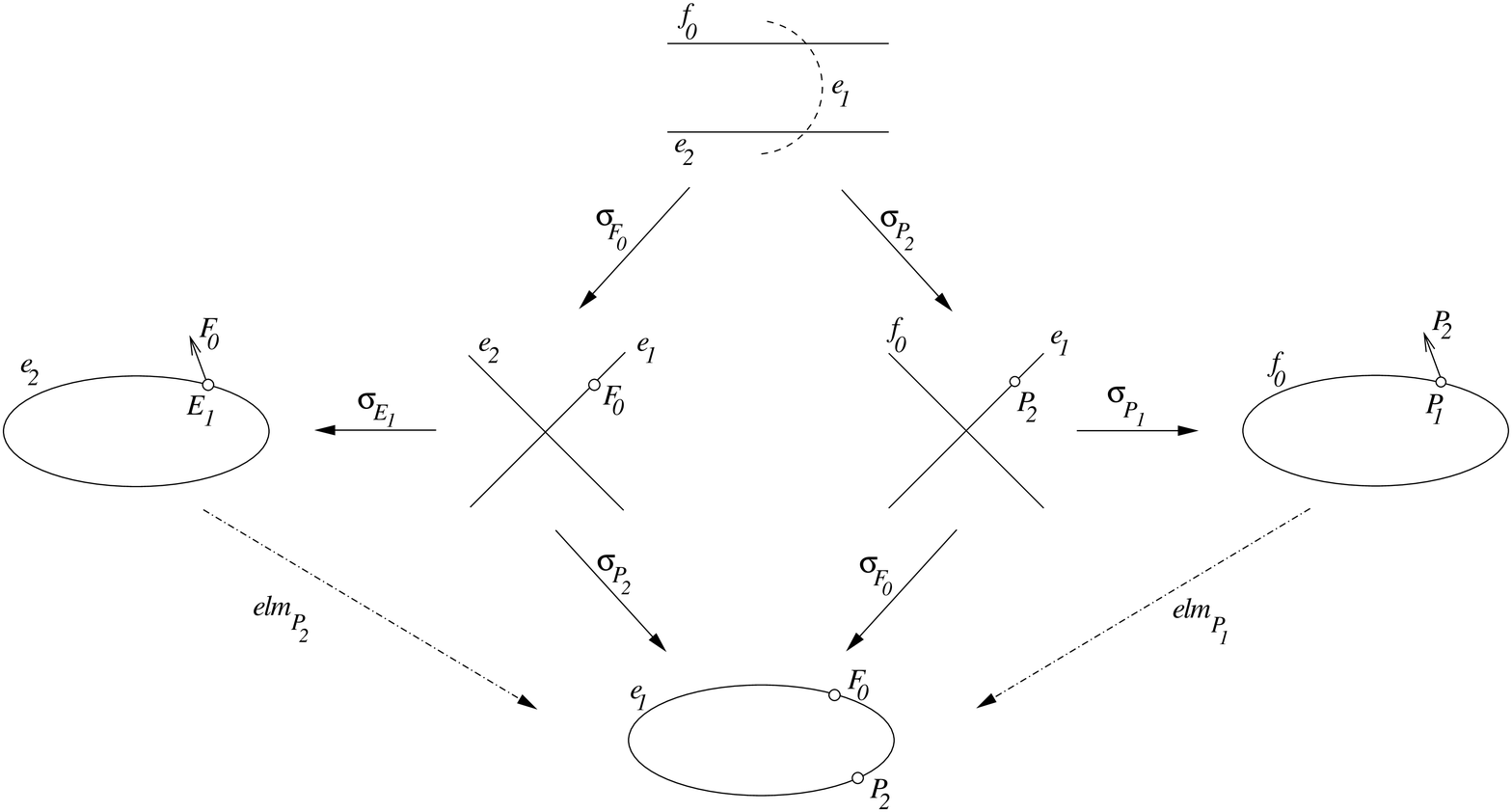}
 }
 \centerline {Figure 11}
 \medskip\noindent
Therefore $S'_0 = elm_{P_1}(S_0)$ and $S''_0 = elm_{P_2}(elm_{P_1}(S_0))$.

\medskip
From the two examples above it is clear that, 
given a main fibre $f_n(A)=~ f_0 \cup e_1 \cup e_2
\cup \dots \cup e_n$, we obtain $n+1$ geometrically ruled surfaces by contracting $n$
components of it. Moreover, if one of them is $S_0$, then  each other is obtained from $S_0$ by
a chain of elementary transformations. \hb
From an analogous procedure,  one can see the general behaviour also in case $f_n(D)$:  the only
possibilities to contract it to a smooth conic are
$con(e_n,e_{n-1}, \dots, e_2,e_1)$ and $con(e_n,e_{n-1}, \dots, e_2,f_0)$. Therefore, the only
geometrically ruled surfaces giving rise to $\wss_0$ are $S_0$ and $elm_{P_1}(S_0)$.

\medskip
The above observations lead to the following result:

\medskip
\proclaim
Theorem 3.3. Let $S$ be a surface ruled by conics and assume that it has a unique
singular fibre $F_n$  embedded of level $n$.  Let $S_0$ and $S'_0$ be two distinct surfaces
geometrically ruled by conics giving rise to $S$ with a minimal number of
blow--ups and contractions. The following facts hold:
\item {-} if $F_n$ is of type $F_n(A)$ then $S'_0 \in \{elm_{P_1}(S_0),
elm_{P_1P_2}(S_0), \dots,  elm_{P_1P_2 \dots P_n}(S_0)\}$, where $P_1 \in S_0$ is a suitable
point and each $P_i$ is a suitable  point, infinitely near to $P_1$  of order $i-1$;
\item {-} if $F_n$ is of type $F_n(D)$ then $S'_0 = elm_{P_1}(S_0)$, where $P_1 \in S_0$ is
a suitable point.
\cvd
\par

\goodbreak

\medskip
\noindent
{\bf Remark 3.4.} Clearly, if $S$ has more than one degenerate fibre, then the previous
theorem can be generalized in a obvious way. For instance, if $S$ has two degenerate fibres
$F_n(A)$ and
$F_m(A)$ and $S_0, S'_0$ are two distinct geometrically ruled
by conics surfaces as in 3.3, then $S'_0$ belongs to the set 
$$
\{elm_\Sigma(S_0) \; | \; \Sigma \subseteq \{P_1, \dots, P_n, Q_1, \dots, Q_m\} \} 
$$
where $P_1, Q_1 \in S_0$ are suitable points, $P_2, \dots, P_n$ are suitable points, infinitely
near to
$P_1$, and 
$Q_2,
\dots, Q_m$ are suitable points, infinitely near to $Q_1$. Clearly, $\Sigma$ has to fulfil the
requirement: if $P_i$ (resp. $Q_i$) $\in \Sigma$ then $P_h$ (resp. $Q_h$) $\in \Sigma$ for all
$h <i$.

\medskip
\noindent
{\bf Definition 3.5.} Let $F^{(1)}, \dots, F^{(p)}$ be the degenerate fibres of $S$ and
let
$l_i$ be the level of $F^{(i)}$, for $i=1, \dots, p$. 
If $\sum_{i=1}^p l_i =L$, we say that $S$ is of {\it level} $L$.

\bigskip
The results concerning one singular fibre (3.3) and two singular fibres (3.4) can be
easily generalized as follows.

\medskip
\noindent
{\bf Remark 3.6.} Let $S$ be a surface, ruled by conics of level $L$ and  
$S_0, S'_0$ be two distinct surfaces, geometrically ruled
by conics  as in 3.3. Then $S'_0$ belongs to the set
$$
 \{elm_\Sigma(S_0) \; | \; \Sigma \; \hbox{is a suitable set of points and } |\Sigma|
\le L \}.
$$
Conversely, note that each surface in the above set is geometrically ruled by conics and
gives rise to $S$ with exactly $L$ blow--ups (followed by contraction). Therefore, the
above set coincides with the following:
$$
\{S_0 \; | \; \hbox{$S_0$ is a g.r.s.  and $S$ can be obtained from it by a
sequence of $L$ blow--ups and contractions} \; \}.
$$

\medskip
\noindent
{\bf Definition 3.7.} We denote the above set of geometrically ruled surfaces of level
$L$ by
$\grc_L(S)$.

\medskip
All the previous constructions can be interpreted in terms of projections. 
More precisely, we show that $S = \Phi_{\wdd - B}(\wss_0) \subset \Bbb P^M$ can be
obtained from  $S_0 = \Phi_D(\Bbb P({\cal F}))\subset \Bbb P^N$ by a linear projection
from a suitable centre. In the following pictures, $\Phi$ stands for $\Phi_{\wdd - B}$.

\medskip 
 \noindent
{\bf Example 3.8.}    A fibre $F_2(A) \subset S$ can be obtained  as the projection of $f_0$ from two points,
either distinct or coincident. Let us consider the first case: $P_1, P_2 \in f_0$ and
$P_1 \ne  P_2$. The following picture illustrates that $\pi_{P_1 P_2}$ factors trough
the projection $\pi_{P_1}$ (giving a fibre of type $F_1$).

 \medskip    
   \centerline{
 \epsfxsize=10cm 
 \epsfbox{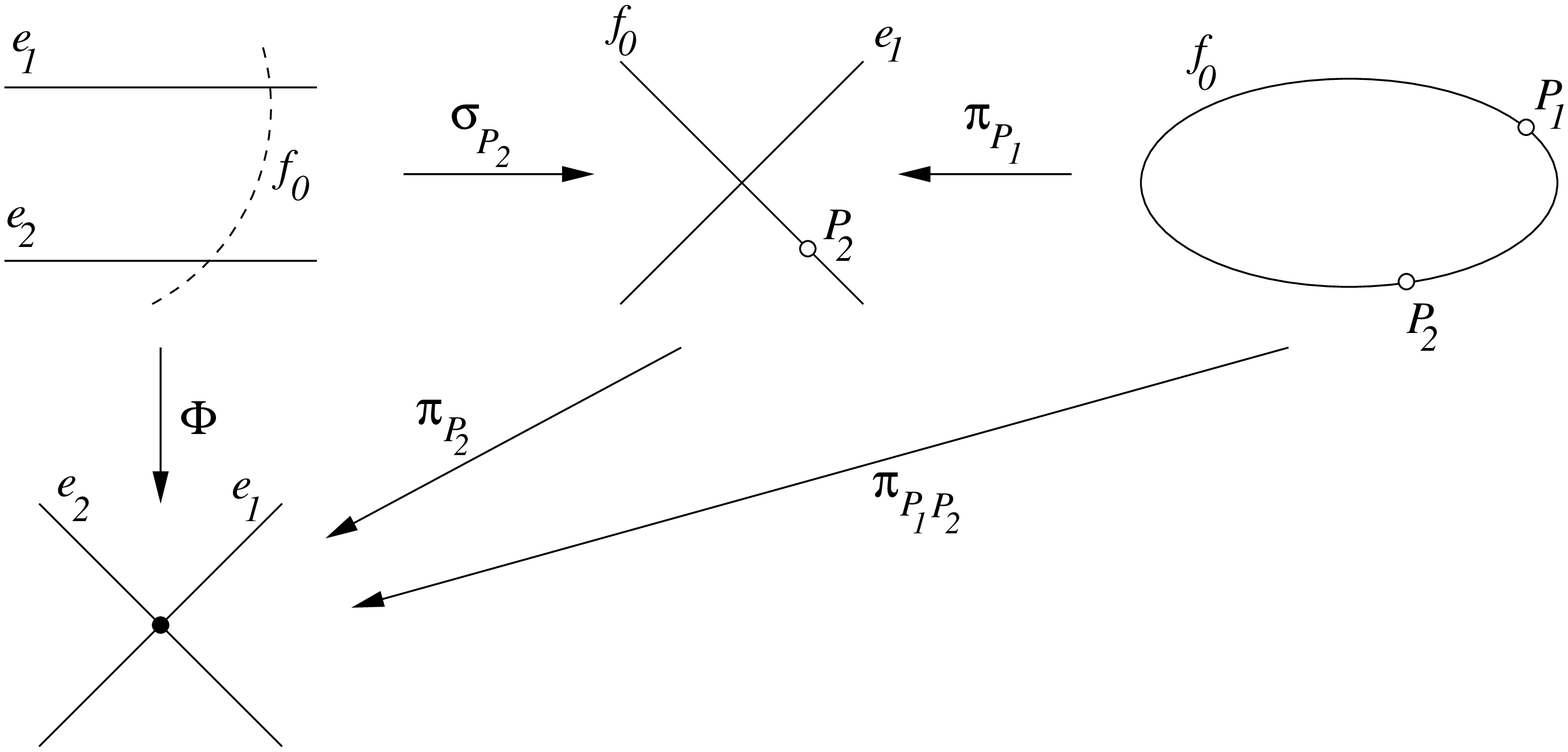}
 }
 \centerline {Figure 12}
 \medskip

\noindent
It is clear that $\pi_{P_1 P_2}$ contracts $f_0$ to a point. This is due to the fact that $\pi_{P_1 P_2}$ is
the projection centered in the line $P_1P_2$ lying on the plane spanned by the conic $f_0$. \hb
The same argument runs in the case $P_1 \in f_0$ and $P_2$ infinitely near point {\it along a transversal
direction}.

\medskip
Similar construction can be performed in the general case $F_n(A)$ 
as the following  picture shows.
\medskip    
   \centerline{
 \epsfxsize=10cm 
\epsfbox{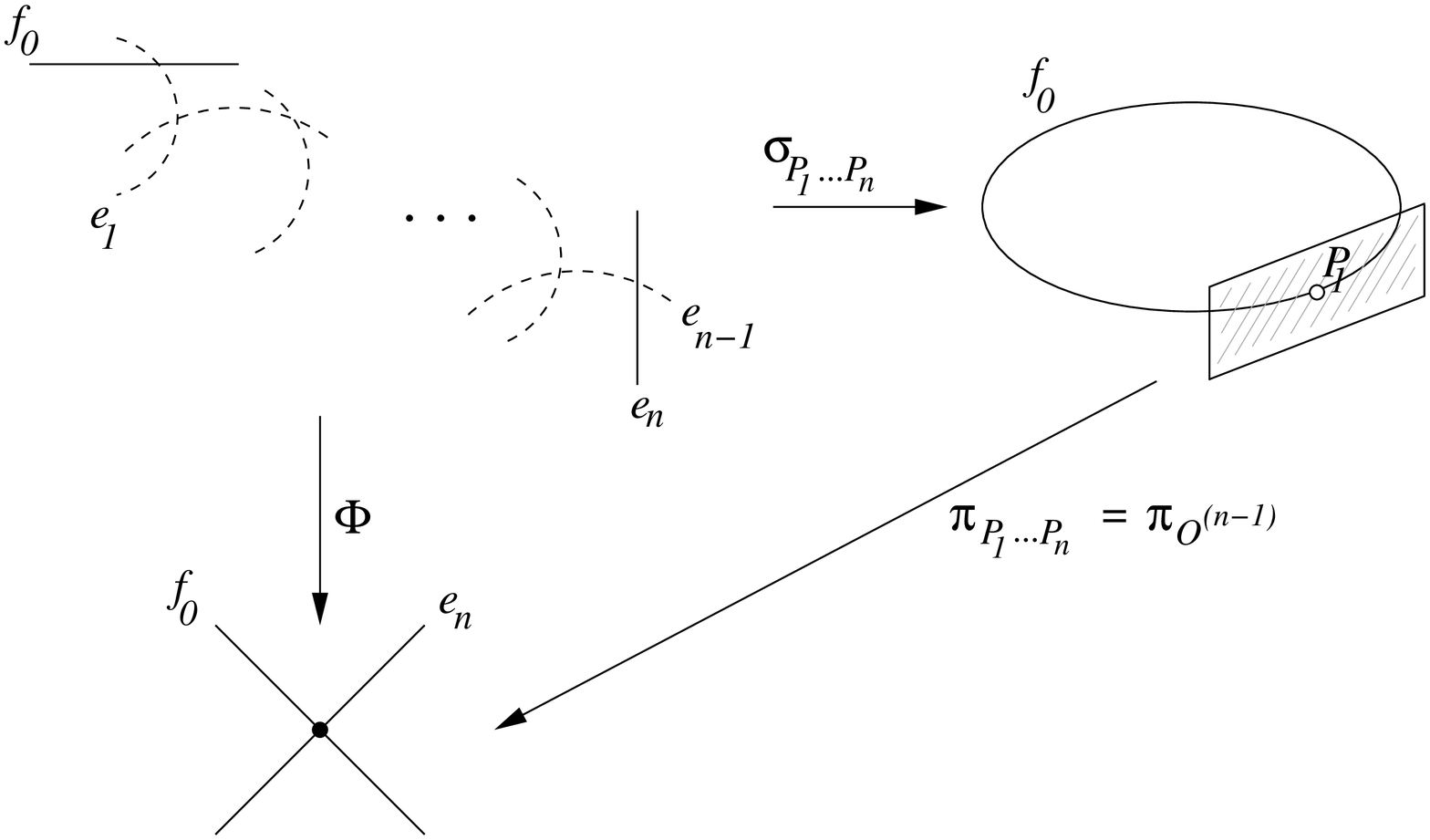}
 }
 \centerline {Figure 13}
 \medskip
 \noindent
and an analogous argument holds for  $F_n(D)$. In these cases, $O^{(n-1)}$ denotes the osculating space of
dimension $n-1$ to a suitable unisecant curve $U \subset S_0$ passing through $P_1$ and
such that $P_2, \dots, P_n$ are points infinitely near to $P_1$ along $U$. In this way,  it is immediate to show the following fact:

\medskip
\noindent
\proclaim 
Proposition 3.9. Let  $S \subset \Bbb P^M$ be a surface ruled by conics of level $L$
and $S_0 \subset \Bbb P^N$ be a surface  in $\grc_L(S)$. Then there exists a
projection $\harrowlength=25pt \pi: \Bbb P^N
\mapright \Bbb P^M$ such that $S = \pi(S_0)$ and $\deg(S) = \deg(S_0) - L$.
\cvd

\bigskip
\centerline {\bf Acknowledgements } 
\smallskip
The referee of an earlier version of [3] pointed out that a
certain argument we gave in that paper was insufficiently clear. His
observations led us to begin this further, more complete study: we
are grateful to that referee for his accuracy.

We are deeply indebted to Peter Newstead and M. S. Narasimhan for their warm support and helpful suggestions; we are   also very grateful to Valentina Beorchia and Dario Portelli  for many interesting discussions.

\bigskip
\centerline{\bf References}

\medskip

\item{[1]}  Badescu, L. (2001). {\it Algebraic Surfaces}. Universitext. New York:  Springer--Verlag.

\medskip

\item{[2]}   Brundu,  M., Sacchiero, G.  (2003). On rational surfaces ruled by conics. {\it
Comm. Algebra } 31, 8: 3631--3652.

\medskip

\item{[3]}  Brundu,  M., Sacchiero, G.   Stratification of the moduli space of
four--gonal curves. {\it Submitted}.

\medskip

\item{[4]}   Chan, D.,  Ingalls, C. (2012). Conic bundles and Clifford algebras.  {\it Contemporary Math.}  562: 53-76.

\medskip

\item{[5]}  Friedman, R.  (1998). {\it Algebraic Surfaces and Holomorphic Vector Bundles}.
Universitext. New York:  Springer--Verlag.

\medskip

\item{[6]} Hartshorne, R. (1977). {\it Algebraic Geometry}. Graduate Texts in Mathematics 52. New York:  Springer--Verlag. 

\medskip

\item{[7]}  Iskovskikh, V.A. (1996).  On a rationality criterion for conic bundles.   {\it Sb. Math.} 187, 7:  1021--1038.

\medskip

 \item{[8]}  Mori, S., Prokhorov,  Y.  (2008). On $\Bbb Q$--conic bundles.  {\it Publ. Res. Inst. Math. Sci.}  44,  2:  315--369.

\medskip

\item{[9]}  Sarkisov, V.G. (1983). On conic bundle structures. {\it Math. USSR-Izv.}  20, 2: 355--390. 

\medskip

\item{[10]} Sarkisov, V.G. (1981). Birational automorphisms of conic bundles.   {\it Math. USSR-Izv.} 17, 4:  177--202.

\medskip

\item{[11]}  Serre, J.P. (2000). {\it Local Algebra}. New York:  Springer--Verlag.

\bye